\documentclass{elsarticle}

\usepackage[margin=1in]{geometry}

\usepackage{graphicx}
\usepackage{pstool}
\usepackage{wrapfig}
\usepackage{caption}
\usepackage{subcaption}
\usepackage[section]{placeins} % Defines \FloatBarrier which does not allow
\usepackage{fancyhdr} % Required for custom headers
\usepackage{lastpage} % Required to determine the last page for the footer
\usepackage{extramarks} % Required for headers and footers

\usepackage{xcolor} % Colors
\usepackage{enumerate} % Additional options for \enumerate{}
\usepackage{paralist} % Additional options for lists; inline lists
\usepackage{amsmath, amsthm, amssymb, mathtools}
\usepackage{mathabx, pifont, stmaryrd} % More math symbols, fonts
\usepackage[explicit]{titlesec} % Control over section/subsection styles
\usepackage{etoolbox} % Boolean variables and more
\usepackage{relsize}
\usepackage{bibentry} % Bibliography entries anywhere in doc
\makeatletter\let\saved@bibitem\@bibitem\makeatother % Do not remove if bibentry
\usepackage[colorlinks, bookmarksopen, bookmarksnumbered,
            citecolor=red,urlcolor=red]{hyperref} % Hyperlinks
\makeatletter\let\@bibitem\saved@bibitem\makeatother % Do not remove if bibentry
\usepackage{algorithm}
\usepackage{algorithmic}

\usepackage{bm} % required for: \bm{}
\usepackage{amsfonts} % required for: \mathbb{}, \mathcal{}, ...

\DeclareMathOperator*{\argmin}{arg\,min}

\newcommand{\ds}[1]{\ensuremath{\displaystyle{#1}}}
\newcommand{\func}[3]{\ensuremath{#1 : #2 \rightarrow #3}}

\newcommand{\norm}[1]{\ensuremath{\left\| #1 \right\|}}

\newcommand{\suchthat}{\mathrel{}\middle|\mathrel{}}

\newcommand{\pder}[2]{\ensuremath{\frac{\partial #1}{\partial #2}}}

\newcommand{\Ecal}{\ensuremath{\mathcal{E}}}

\newcommand{\Gcal}{\ensuremath{\mathcal{G}}}
\newcommand{\Hcal}{\ensuremath{\mathcal{H}}}

\newcommand{\Lcal}{\ensuremath{\mathcal{L}}}

\newcommand{\Pcal}{\ensuremath{\mathcal{P}}}
\newcommand{\Qcal}{\ensuremath{\mathcal{Q}}}

\newcommand{\Vcal}{\ensuremath{\mathcal{V}}}
\newcommand{\Wcal}{\ensuremath{\mathcal{W}}}

\newcommand{\Rbb}{\ensuremath{\mathbb{R} }}

\newcommand\Abm{{\ensuremath{\bm{A}}}}

\newcommand\Fbm{{\ensuremath{\bm{F}}}}

\newcommand\Mbm{{\ensuremath{\bm{M}}}}

\newcommand\Rbm{{\ensuremath{\bm{R}}}}

\newcommand\bbm{{\ensuremath{\bm{b}}}}
\newcommand\cbm{{\ensuremath{\bm{c}}}}

\newcommand\fbm{{\ensuremath{\bm{f}}}}

\newcommand\kbm{{\ensuremath{\bm{k}}}}

\newcommand\mbm{{\ensuremath{\bm{m}}}}

\newcommand\rbm{{\ensuremath{\bm{r}}}}

\newcommand\ubm{{\ensuremath{\bm{u}}}}

\newcommand\wbm{{\ensuremath{\bm{w}}}}
\newcommand\xbm{{\ensuremath{\bm{x}}}}
\newcommand\ybm{{\ensuremath{\bm{y}}}}
\newcommand\zbm{{\ensuremath{\bm{z}}}}

\newcommand\lambdabold{{\ensuremath{\boldsymbol{\lambda}}}}

\newcommand\xibold{{\ensuremath{\boldsymbol{\xi}}}}
\newcommand\zetabold{{\ensuremath{\boldsymbol{\zeta}}}}

\newcommand\zerobold{\ensuremath{\mathbf{0}}}

\usepackage{tikz}
\usepackage{pgfplots}
\usepackage{pgfplotstable, filecontents, booktabs}
\pgfplotsset{compat=1.9}

\usetikzlibrary{pgfplots.groupplots}
\usepgfplotslibrary{fillbetween}
\usetikzlibrary{calc,fit,matrix,arrows,automata,positioning,shapes}
\usetikzlibrary{arrows.meta}

\pgfplotsset{select coords between index/.style 2 args={
    x filter/.code={
        \ifnum\coordindex<#1\fi
        \ifnum\coordindex>#2\fi
    }
}}

\tikzset{
 invisible/.style={opacity=0},
 visible on/.style={alt={#1{}{invisible}}},
 alt/.code args={<#1>#2#3}{%
   \alt<#1>{\pgfkeysalso{#2}}{\pgfkeysalso{#3}}
 },
}

\newcommand{\logLogSlopeTriangle}[5]
{
 \pgfplotsextra
 {
  \pgfkeysgetvalue{/pgfplots/xmin}{\xmin}
  \pgfkeysgetvalue{/pgfplots/xmax}{\xmax}
  \pgfkeysgetvalue{/pgfplots/ymin}{\ymin}
  \pgfkeysgetvalue{/pgfplots/ymax}{\ymax}

  \pgfmathsetmacro{\xArel}{#1}
  \pgfmathsetmacro{\yArel}{#3}
  \pgfmathsetmacro{\xBrel}{#1-#2}
  \pgfmathsetmacro{\yBrel}{\yArel}
  \pgfmathsetmacro{\xCrel}{\xArel}
  \pgfmathsetmacro{\lnxB}{\xmin*(1-(#1-#2))+\xmax*(#1-#2)} % in [xmin,xmax].
  \pgfmathsetmacro{\lnxA}{\xmin*(1-#1)+\xmax*#1} % in [xmin,xmax].
  \pgfmathsetmacro{\lnyA}{\ymin*(1-#3)+\ymax*#3} % in [ymin,ymax].
  \pgfmathsetmacro{\lnyC}{\lnyA+#4*(\lnxA-\lnxB)}
  \pgfmathsetmacro{\yCrel}{\lnyC-\ymin)/(\ymax-\ymin)} % THE IMPROVED EXPRESSION WITHOUT 'DIMENSION TOO LARGE' ERROR.

  \coordinate (A) at (rel axis cs:\xArel,\yArel);
  \coordinate (B) at (rel axis cs:\xBrel,\yBrel);
  \coordinate (C) at (rel axis cs:\xCrel,\yCrel);

  \draw[#5]   (A)-- node[pos=0.5,anchor=north,scale=0.8] {1}
              (B)--
              (C)-- node[pos=0.5,anchor=west,scale=0.8] {#4}
              cycle;
 }
}

\usepackage{setspace}
\newbool{fastcompile}
\setbool{fastcompile}{false}

\newcommand{\blacksolidline}{\raisebox{2pt}{\tikz{\draw[-,black,line width = 1.0pt](0,0) -- (5mm,0);}}}

\newcommand{\reddottedline}{\raisebox{2pt}{\tikz{\draw[red,dotted,line width = 1pt,line cap=round, dash pattern=on 0pt off 2\pgflinewidth](0,0) -- (5mm,0);}}}
\newcommand{\reddashedline}{\raisebox{2pt}{\tikz{\draw[-,red,dashed,line width = 0.6pt](0,0) -- (5mm,0);}}}

\newcommand{\bluesolidthinline}{\raisebox{2pt}{\tikz{\draw[-,blue,line width = 0.5pt](0,0) -- (5mm,0);}}}
\newcommand{\redsolidthinline}{\raisebox{2pt}{\tikz{\draw[-,red,line width = 0.5pt](0,0) -- (5mm,0);}}}
\newcommand{\greensolidthinline}{\raisebox{2pt}{\tikz{\draw[-,green,line width = 0.5pt](0,0) -- (5mm,0);}}}
\newcommand{\blacksolidthinline}{\raisebox{2pt}{\tikz{\draw[-,black,line width = 0.5pt](0,0) -- (5mm,0);}}}

\newcommand{\mapped}[1]{{#1}_X}

\newcommand{\defgrad}{G}
\newcommand{\defjac}{g}
\newcommand{\dom}{\Omega}
\newcommand{\refdom}{\Omega_0}
\newcommand{\stvc}{U}
\newcommand{\flux}{F}

\newcommand{\numflux}{\Hcal}

\newcommand{\msh}{\mathrm{msh}}

\newcommand{\adv}{\mathrm{adv}}
\newcommand{\up}{\mathrm{up}}

\begin{document}

\title{Implicit shock tracking for unsteady flows by the method of lines}

\author[rvt1]{A.~Shi\fnref{fn1}\corref{cor1}}
\ead{andrewshi94@berkeley.edu}

\author[rvt1,rvt2]{P.-O.~Persson\fnref{fn2}}
\ead{persson@berkeley.edu}

\author[rvt3]{M.~J.~Zahr\fnref{fn3}}
\ead{mzahr@nd.edu}

\address[rvt1]{Department of Mathematics, University of California, Berkeley,
                    Berkeley, CA 94720, United States}
\address[rvt2]{Mathematics Group, Lawrence Berkeley National Laboratory,
               1 Cyclotron Road, Berkeley, CA 94720, United States}
\address[rvt3]{Department of Aerospace and Mechanical Engineering, University
               of Notre Dame, Notre Dame, IN 46556, United States}
\cortext[cor1]{Corresponding author}
\fntext[fn1]{Graduate student, Department of Mathematics, University of
             California, Berkeley}
\fntext[fn2]{Professor, Department of Mathematics, University of
             California, Berkeley}
\fntext[fn3]{Assistant Professor, Department of Aerospace and Mechanical
             Engineering, University of Notre Dame}


\begin{abstract}
A recently developed high-order implicit shock tracking (HOIST) framework for resolving
discontinuous solutions of inviscid, steady conservation laws
\cite{zahr2018optimization, zahr2020implicit} is extended to
the \textit{unsteady} case. Central to the framework is an optimization
problem which simultaneously computes a  discontinuity-aligned mesh and
the corresponding high-order approximation to the flow, which provides
nonlinear stabilization and a high-order approximation to the solution.
This work extends the implicit shock tracking framework to the case of unsteady conservation laws
using a method of lines discretization via a diagonally implicit Runge-Kutta
method by ``solving a steady problem at each timestep". 
We formulate and solve an optimization problem that produces a
feature-aligned mesh and solution at each Runge-Kutta stage of each timestep,
and advance this solution in time by standard Runge-Kutta update formulas.
A Rankine-Hugoniot based prediction of the shock location together with
a high-order, untangling mesh smoothing procedure provides a high-quality
initial guess for the optimization problem at each time, which results in
rapid convergence of the sequential quadratic programing (SQP)
optimization solver. This method is shown to deliver highly accurate
solutions on coarse, high-order discretizations without nonlinear
stabilization and recover the design accuracy of the Runge-Kutta
scheme. We demonstrate this framework on a series of inviscid, unsteady
conservation laws in both one- and two- dimensions. We also 
verify that our method is able to recover the design order of accuracy
of our time integrator in the presence of a strong discontinuity.
\end{abstract}

\maketitle

\section{Introduction}
\label{sec:intro}

It is widely believed that higher fidelity is required for problems with
propagating waves, turbulent fluid flow, nonlinear interactions, and multiple
scales \cite{wang2013high}. This has resulted in a significant interest in
high-order accurate methods, such as discontinuous Galerkin (DG) methods
\cite{cockburn01rkdg,hesthaven2007nodal}, which have the potential to produce
accurate solutions on coarse meshes. Among the most significant challenges
associated with high-order methods is their sensitivity to under-resolved
features, in particular for nonlinear problems where the spurious
oscillations often cause a breakdown of the numerical solvers. This is
exacerbated for problems with shocks where the low dissipation associated
with high-order methods is insufficient to stabilize the solution. Since
shocks are present in many important problems in fields such as aerospace,
astrophysics, and combustion, this poses a fundamental barrier to widespread
adoption of these methods. 

Several approaches have been proposed to stabilize shocks, 
most of which are based on \emph{shock capturing}, where the numerical 
discretization accounts for discontinuities independently of the
computational grid.  One simple method is to use a sensor that identifies
the mesh elements  in the shock region and reduce the degree of the
approximating polynomial \cite{BauOden,burbeau01limiter}. A more
sophisticated approach includes limiting,  such as the weighted
essentially non-oscillatory (WENO) schemes \cite{eno1,weno1,weno2},
which perform a high-order reconstruction near discontinuities, but
require a large computational stencil  which offsets the attractive
compactness properties of DG methods. For high-order methods, artificial
viscosity has also proven to be competitive, since it can smoothly resolve
the jumps in the solution without introducing additional discontinuities
between the elements \cite{persson06shock}. A recent comparative study of
artificial viscosity models \cite{yu2020study} discusses their relative merits,
but notes they all suffer from a relatively strong dependency on a large 
number of empirical parameters which must be tuned.
The main problem with all these approaches
is they lead to globally first-order accurate schemes. This can be remedied
by local mesh refinement around the shock ($h$-adaptivity)
\cite{dervieux03adaptation}, although the anisotropic high-order mesh
adaptation is challenging and requires highly refined elements near
the shock. This issue is further complicated by transient problems where shocks 
and other local features propagate throughout the domain, requiring online
adaptivity to be computationally feasible. Even with aggressive $h$- and $p$-adaptivity, 
simulations of complex, shock-dominated, unsteady flows are extremely challenging and expensive.

%1.3) Previous Shock Tracking Approaches & Explicit Shock Tracking
An alternative approach is \emph{shock tracking} or \emph{shock fitting},
where the computational mesh is moved such that its faces are aligned with
the discontinuities in the solution. This is natural in the setting of a DG
method since the numerical scheme already incorporates jumps between the
elements and the approximate Riemann solvers employed on the element faces
handle the discontinuities correctly. However, it is a difficult meshing
problem since it essentially requires generating a fitted mesh to the
(unknown) shock surface. Many previous approaches employ specialized
formulations and solvers which are dimension dependent and do not easily
generalize \cite{harten1983self, glimm2003conservative, bell1982fully}
and/or are limited to relatively simple problems
\cite{shubin1981steady, shubin1982steady, rosendale1994floating}.
In addition, early approaches to shock fitting have been applied to
low-order schemes where the relative advantage over shock capturing is
smaller than for high-order methods
\cite{trepanier1996conservative, baines2002multidimensional}. 
One particular class of methods of note is \textit{explicit} shock tracking,
which is surveyed in \cite{moretti2002thirty, salas2009shock}.
These strategies largely consist of explicitly identifying the
shock and using the Rankine-Hugoniot conditions
to compute its motion and states upstream and downstream of the shock.
More recent developments in explicit shock tracking \cite{rawat2010high}
use more sophisticated methods to compute shock velocities and 
discretize the flow equations, but ultimately still require a specialized
strategy  to explicitly track the shock separately from the remainder of the
flow. These methods are not easily applicable to 
discontinuities whose topologies not known \textit{a priori}. 
While interest in shock tracking/fitting methods has seen
somewhat of a resurgence in recent years 
\cite{Ciallella2020ExtrapolatedST, Bonfiglioli2015shkfit, Geisenhofer2020Extended, DAquila2020shkfit}, 
shock tracking is largely not used in practical CFD today. 

%1.4) 2018/2020 Shock Tracking
In \cite{zahr2018optimization, zahr2020implicit}, we introduced a novel
approach to shock tracking for steady conservation laws that does not require
explicitly generating a mesh of the unknown discontinuity surface. Rather,
the conservation law is discretized on a mesh without knowledge of the
discontinuity surface and an optimization problem is formulated such that
its solution is the pair $(\ubm,\xbm)$, where $\xbm$ is the position of the
mesh nodes that cause element faces to align with discontinuities in the
flow and $\ubm$ is the solution of discretized conservation law on the mesh
defined by $\xbm$. That is, discontinuity tracking is implicitly achieved
through the solution of an optimization problem and will be referred to as
\textit{implicit shock tracking}. While this approach works with
any discretization that allows for inter-element discontinuities,
we focus on high-order DG methods due to the high degree of accuracy
attainable on coarse meshes, proper treatment of discontinuities with
approximate Riemann solvers, and the ability to use curved elements to
track discontinuities with curvature.
The implicit tracking optimization problem proposed in
\cite{zahr2020implicit} minimizes the violation of the DG
residual in an enriched test space while enforcing that the
standard DG (same test and trial space) equation is satisfied.
This objective function is a surrogate for violation of the
infinite-dimensional weak formulation of the conservation law,
which endows the method with $r$-adaptive behavior: it promotes
alignment of the mesh with discontinuities and adjusts nodes in
smooth regions to improve approximation of the conservation law.
The optimization problem is solved using a sequential quadratic
programming method with a Levenberg-Marquardt Hessian approximation
that simultaneously converges the mesh and solution
to their optimal values, which never requires the fully converged
DG solution on a non-aligned mesh and does not require nonlinear
stabilization. The combination of implicit tracking with a DG
discretization leads to a high-order accurate numerical method
that has been shown to provide accurate approximations to high-speed inert
\cite{zahr2018optimization, zahr2020implicit} and
reacting flows \cite{zahr2020react}. 

%1.5) This MOL Paper
In this work, we further extend the framework developed in
\cite{zahr2018optimization, zahr2020implicit} for steady conservation laws
(which can be applied to space-time formulations of unsteady conservation laws)
to a method of lines discretization approach for unsteady problems. 
While space-time methods are attractive for a number of reasons, the method of
lines approach tends to be more practical for complex problems when applicable,
in large part because computations are only required on a $d$-dimensional
mesh instead of a time-coupled $(d+1)$-dimensional space-time mesh.
The key ingredients of the method of lines approach are:
\begin{inparaenum}[1)]
\item an Arbitrary Lagrangian-Eulerian formulation of the conservation law
      to handle the deforming mesh (which deforms to track the shock through
      the domain),
\item semi-discretization with DG to obtain a system of ordinary differential
      equations,
\item high-order temporal discretization with a diagonally implicit Runge-Kutta (DIRK) method, and
\item implicit shock tracking at each time step following the
      approach in \cite{zahr2020implicit}.
\end{inparaenum}
We utilize a Rankine-Hugoniot-based procedure to predict the shock location 
at future times combined with a high-order mesh smoothing procedure 
with untangling capabilities to construct high quality initial guesses
for the mesh and solution for the optimization problem at each time step.
These initial guesses enhance the robustness of the method and 
significantly accelerate the performance of the optimization solver,
to achieve rapid convergence.

%1.6) Related NRL Work
To our knowledge, the only other approach to implicit shock tracking is
the Moving Discontinuous Galerkin Method with Interface Condition Enforcement
(MDG-ICE), proposed in \cite{corrigan2019moving, corrigan2019convergence},
where the authors enforce a DG discretization with unconventional 
numerical fluxes and the Rankine-Hugoniot interface conditions 
in a minimum-residual sense. In their approach, the interface condition
is enforced along all faces of the mesh, separately from conservation law.
Interestingly, enforcement of the interface condition circumvented 
traditional stability requirements for the DG numerical fluxes, 
allowing them to solely rely on fluxes interior
to an element. Their method was shown to successfully track even complex
discontinuity surfaces and provide high-order approximations to the
conservation law on traditionally coarse, high-order meshes.
In \cite{corrigan2019unsteady}, they extended MDG-ICE to solve
unsteady problems by space-time slab marching with a simplicial grid extrusion
strategy. More recently, they extended MDG-ICE, originally developed for inviscid
conservation laws, to viscous conservation laws \cite{kercher_moving_2021} and
reformulated it as a least-squares discontinuous Galerkin method
\cite{kercher_least-squares_2020}, which endows the method with
super-optimal convergence properties.

%1.7) Organization of this Paper
The remainder of this paper is organized as follows. Section~\ref{sec:disc}
introduces the governing system of inviscid unsteady conservation laws,
its spatial discretization using a DG method, and its temporal discretization 
using a diagonally implicit Runge-Kutta (DIRK) method.
Section~\ref{sec:optim} presents the error-based objective function 
and the constrained optimization framework. 
Section~\ref{sec:practical} discusses practical details required for the 
proposed tracking framework such as initialization of the SQP solver.
Finally, Section~\ref{sec:num-exp} presents a number of numerical
experiments that demonstrate the method on a variety of unsteady 
flows using coarse, high-order meshes. We also demonstrate
high-order temporal convergence of the method along with the
rapid convergence of the SQP solver.

\section{Governing equations and high-order numerical discretization}
\label{sec:disc}

Consider a general system of $m$ inviscid conservation laws, defined on the
fixed domain $\dom \subset \Rbb^d$ and subject to appropriate boundary
conditions,
\begin{equation} \label{eqn:claw-phys}
 \pder{U}{t} + \nabla\cdot \flux(\stvc) = 0
                \quad \text{in}~~\dom\times[0,T]
\end{equation}
where $\func{\stvc}{\dom\times[0,T]}{\Rbb^m}$ is the solution of the system of
conservation laws, $\func{\flux}{\Rbb^m}{\Rbb^{m\times d}}$ is the flux
function,
$\ds{\nabla \coloneqq (\partial_{x_1},\dots,\partial_{x_d})}$
is the gradient operator in the physical domain such that
$\nabla W(x,t) =
 \begin{bmatrix} \partial_{x_1} W(x,t) & \cdots &
                 \partial_{x_d} W(x,t)
 \end{bmatrix} \in \Rbb^{N \times d}$
for any $\func{W}{\Omega\times[0,T]}{\Rbb^N}$ and $x\in\Omega$, $t\in[0,T]$,
and the boundary of the domain $\partial\Omega$ has outward unit normal
$\func{n}{\partial\dom}{\Rbb^d}$. The conservation law
in (\ref{eqn:claw-phys}) is supplemented with the initial
condition $U(x, 0)=\bar{U}(x)$ for all $x\in\Omega$, where
$\func{\bar{U}}{\Omega}{\Rbb^m}$.
In general, the solution $\stvc(x)$ may contain discontinuities,
in which case, the conservation law
(\ref{eqn:claw-phys}) holds away from the discontinuities
and the Rankine-Hugoniot conditions \cite{majda2012compressible}
hold at discontinuities.

Building on our previous work \cite{zahr2018optimization, zahr2020implicit}, 
we will construct a high-order numerical method that tracks
discontinuities with the computational grid as they evolve
through the domain, which places three requirements on the
discretization:
\begin{inparaenum}[1)]
 \item a high-order, stable, and convergent discretization
       of the conservation law in (\ref{eqn:claw-phys}),
 \item employs a solution basis that supports discontinuities between
       computational cells or elements, and
 \item allows for deformation of the computational domain.
\end{inparaenum}
As such, our method is based on a standard high-order DG-DIRK discretization
of an Arbitrary Lagrangian-Eulerian (ALE) formulation of the governing
equations. We insist on a high-order discretization given their proven
ability \cite{zahr2018optimization, corrigan2019moving} to deliver accurate
solutions on coarse discretizations provided discontinuities are tracked.

The remainder of this section will detail the discretization of the
conservation law (\ref{eqn:claw-phys}) using DG such that it reduces
to the semi-discrete form
\begin{equation} \label{eqn:claw-semidisc}
\rbm(\dot\ubm, \ubm, \dot\xbm, \xbm) = \zerobold
\end{equation}
where $\rbm: \Rbb^{N_\ubm} \times \Rbb^{N_\ubm} \times \Rbb^{N_\xbm} \times \Rbb^{N_\xbm} \rightarrow \Rbb^{N_\ubm}$,
$\ubm$ is the semi-discrete representation of the conservation law state $U$,
and $\xbm$ is the semi-discrete representation of the conservation law domain
$\Omega$ (nodal coordinates of mesh nodes). %Usually, $\xbm$ is known
%analytically (e.g., forced motion) or governed by a dynamical system
%(e.g., fluid-structure interaction); however, in the implicit shock
%tracking setting, it is unknown and determined as the solution of
%the central optimization problem.
We will then apply a high-order temporal discretization to 
(\ref{eqn:claw-semidisc}) by a diagonally implicit Runge-Kutta method
to yield a complete discretization of (\ref{eqn:claw-phys}). 
The same discretization process can be used to yield a semi-discretization 
from an enriched \textit{test} space and corresponding temporal
discretization, which will be used in the definition of the proposed
objective function.

\subsection{Arbitrary Lagrangian-Eulerian formulation of conservation laws}
\label{sec:disc:transf}
\begin{figure}
  \centering
  \input{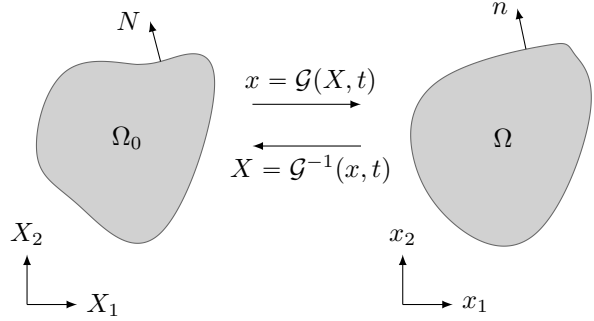}
  \caption{Mapping between reference and physical domains.}
  \label{fig:dom-map}
\end{figure}

We use an ALE formulation of the governing equations to account
for the time-dependent domain deformations required to track
discontinuities as they evolve. To this end, we introduce a
time-dependent domain mapping (Fig.~\ref{fig:dom-map})
\begin{equation} \label{eqn:dom-map}
 \func{\Gcal}{\Omega_0\times[0,T]}{\Omega}; \qquad
 \Gcal : (X,t) \mapsto \Gcal(X,t),
\end{equation}
where $\Omega_0 \subset \Rbb^d$ is a fixed reference domain,
$T$ is the final time, and at each time $t\in[0,T]$,
$\func{\Gcal(\,\cdot\,,t)}{\Omega_0}{\Omega}$ is a
diffeomorphism. We note that the domain $\Omega$ is fixed,
i.e., $\Omega$ occupies the same region of $\Rbb^d$ at any
time $t\in[0,T]$; the time-dependent diffeomorphism is
introduced as an integral part of the proposed numerical
method to track discontinuities as they evolve.
Under the domain mapping (\ref{eqn:dom-map}), 
the conservation law (\ref{eqn:claw-phys}) becomes
\begin{equation} \label{eqn:claw-phys-transfdom}
 \pder{U}{t} + \nabla \cdot F(U) = 0 \quad \text{in}~~\Gcal(\Omega_0, t)
\end{equation}
Following the approach in \cite{persson2009dgdeform}, the conservation
law on the physical domain $\Omega$ is transformed to a conservation law
on the reference domain $\Omega_0$
\begin{equation} \label{eqn:claw-ref}
 \pder{U_X}{t} + \nabla_X \cdot F_X(U_X; G, v) = 0 \quad \text{in}~~\Omega_0
\end{equation}
where $\func{\stvc_X}{\dom_0\times[0,T]}{\Rbb^m}$ is the solution of the
transformed conservation law,
$\func{\flux_X}{\Rbb^m\times\Rbb^{d\times d}\times\Rbb^d}{\Rbb^{m\times d}}$
is the transformed flux function,
$\ds{\nabla_X \coloneqq (\partial_{X_1},\dots,\partial_{X_d})}$
is the gradient operator on the reference domain, and
the deformation gradient
$\func{G}{\Omega_0\times[0,T]}{\Rbb^{d\times d}}$,
mapping Jacobian
$\func{g}{\Omega_0\times[0,T]}{\Rbb}$, and
mapping velocity
$\func{v}{\Omega_0\times[0,T]}{\Rbb^d}$ are defined as
\begin{equation}
 G = \nabla_X\Gcal, \qquad g = \det G, \qquad v = \pder{\Gcal}{t}.
\end{equation}
The transformed and physical
solutions are related, for any $X\in\Omega_0$ and $t\in[0,T]$, as
\begin{equation} \label{eqn:transf-stvc}
 U_X(X,t) = g(X,t) U(\Gcal(X,t),t)
\end{equation}
and the transformed flux is defined as
\begin{equation} \label{eqn:transf-flux}
 F_X : (W_X; \Theta, \xi) \mapsto
 [(\det\Theta) F((\det\Theta)^{-1} W_X)- W_X \otimes \xi]\Theta^{-T}.
\end{equation}
The unit normals in the reference and physical domain are related by
\begin{equation} \label{eqn:normal-transf}
n = \frac{gG^{-T}N}{\|gG^{-T}N\|}.
\end{equation}
The transformed conservation law is supplemented with the initial
condition $U_X(X,0)=\bar{U}_X(X)$ for all $X\in\Omega_0$, where
$\func{\bar{U}_X}{\Omega_0}{\Rbb^m}$ is
$\bar{U}_X(X) = g(X,0)\bar{U}(\Gcal(X,0))$.
In this work, we take the reference domain to be the physical domain
at time $0$, which implies $g(X, 0) = 1$.

\subsection{Discontinuous Galerkin discretization of
            the transformed conservation law}
\label{sec:disc:dg}
We use a nodal discontinuous Galerkin method
\cite{cockburn01rkdg, hesthaven2007nodal} to discretize
the transformed conservation law (\ref{eqn:claw-ref}).
Let $\Ecal_h$ represent a discretization of the reference domain
$\refdom$ into non-overlapping, potentially curved, computational elements.
The DG construction begins with the elementwise weak form 
of the conservation  law (\ref{eqn:claw-ref}) that results from multiplying 
each equation by a test function $\psi_X$, integrating over a single element 
$K \in \Ecal_h$, and applying the divergence theorem
\begin{equation} \label{eqn:claw-weak-elem}
 \int_K \psi_X \cdot \dot{U}_X \, dV +
 \int_{\partial K} \mapped\psi^+ \cdot \mapped\flux(\mapped\stvc; G, v) N \, dS -
 \int_K F_X(\mapped\stvc; G, v):\mapped\nabla \mapped\psi \, dV = 0,
\end{equation}
where $N$ is the outward normal to the surface $\partial K$ and
$\mapped\psi^+$ denotes the trace of $\psi$ \textit{interior} to
element $K$. To ensure the face integrals are single-valued, we replace
$\mapped\flux(\mapped\stvc; G, v) N$ in the second term with a numerical flux
function $\func{\Hcal_X}{\Rbb^m\times\Rbb^m\times\Rbb^d\times\Rbb^{d\times d}\times\Rbb^d}{\Rbb^m}$ 
associated with the reference inviscid flux $F_X$
\begin{equation} \label{eqn:claw-weak-elem-numflux}
 \int_K \psi_X \cdot \dot{U}_X \, dV +
 \int_{\partial K} \mapped\psi^+ \cdot
   \mapped\numflux(\mapped\stvc^+,\mapped\stvc^-,N;G, v) \, dS -
 \int_K F_X(\mapped\stvc; G, v):\mapped\nabla \mapped\psi \, dV = 0,
\end{equation}
where $\mapped\stvc^+$ ($\mapped\stvc^-$) denotes the interior (exterior)
trace of $\mapped\stvc$ to the element $K$; for points
$X\in\partial K\cap\partial\Omega_0$, $\mapped{\stvc}^-$
is a boundary state constructed to enforce the appropriate
boundary condition.  In this work, we take the numerical flux
to be a smoothed version of the Roe flux (Section \ref{sec:disc:numflux}).
To establish the finite-dimensional (semi-discrete) form of (\ref{eqn:claw-weak-elem-numflux}),
we introduce the DG approximation (trial) space of discontinuous
piecewise polynomials associated with the mesh $\Ecal_h$
\begin{equation}
 \Vcal_h^p = \left\{v \in [L^2(\refdom\times[0,T])]^m \suchthat
         \left.v(\,\cdot\,,t)\right|_K \in [\Pcal_p(K)]^m,
         ~\forall K \in \Ecal_h,~t\in[0,T]\right\}
\end{equation}
where $\Pcal_p(K)$ is the space of polynomial functions of degree at most
$p \geq 1$ on the element $K$, and we take the DG test space to
be $\Vcal_h^{p'}$, where $p'\geq p$.
We also define the space of admissible domain mappings
as the space of continuous piecewise polynomials of degree
$q$ associated with the mesh $\Ecal_h$
\begin{equation} \label{eqn:gfcnsp}
  \Wcal_h = \left\{v \in [C^0(\Omega_0\times[0,T])]^d \suchthat
         \left.v(\,\cdot\,,t)\right|_K \in [\Pcal_q(K)]^d,
         ~\forall K \in \Ecal_h,~t\in[0,T]\right\}.
\end{equation}
We can now establish the finite dimensional form of 
(\ref{eqn:claw-weak-elem-numflux}) and the formal statement of DG as: 
given $\Gcal_h \in \Wcal_h$, find $U_{X, h} \in  \Vcal_h^p$ such that for all
$\psi_{X, h} \in \Vcal_h^{p'}$
\begin{equation} \label{eqn:claw-weak-elem-finitedim}
 \int_K \psi_{X, h}  \cdot \dot{U}_{X, h} \, dV +
 \int_{\partial K} \psi_{X, h}^+ \cdot
   \mapped\numflux(U_{X, h}^+,U_{X, h}^-,N;G_h, v_h) \, dS -
 \int_K F_X(U_{X, h}; G_h, v_h):\mapped\nabla \psi_{X, h} \, dV = 0,
\end{equation}
where the DG residual form $\func{r_h^{p'p}}{\Vcal_h^{p'}\times\Vcal_h^p\times\Wcal_h}{\Rbb}$ is given by
\begin{equation}
 r_h^{p',p} : (\psi_{X, h}, W_{X, h}, \Qcal_h) \mapsto
 \sum_{K\in\Ecal_{h,q}} r_K^{p',p}(\psi_{X, h}, W_{X, h}, \Qcal_h),
\end{equation}
and the elemental DG form
$\func{r_K^{p',p}}{\Vcal_h^{p'}\times\Vcal_h^p\times\Wcal_h}{\Rbb}$
is given by
\begin{equation} \label{eqn:claw-weak-numflux}
\begin{aligned}
 r_K^{p',p} : (\psi_{X, h}, W_{X, h}, \Qcal_h) \mapsto
 &\int_K \psi_{X, h} \cdot \dot{W}_{X, h} \, dV \\
 &+ \int_{\partial K} \psi_{X, h}^+ \cdot
   \mapped\numflux(W_{X, h}^+,W_{X, h}^-,N;\nabla_X\Qcal_h,\dot\Qcal_h) \, dS \\
 &-\int_K \mapped\flux(W_{X, h};\nabla_X\Qcal_h,\dot\Qcal_h):\mapped\nabla \psi_{X, h} \, dV.
\end{aligned}
\end{equation}
Next, we introduce a (nodal) basis over each element
for the test space ($\Vcal_h^p$), trial space ($\Vcal_h^{p'}$), and
mapping space ($\Wcal_h$) to reduce (\ref{eqn:claw-weak-numflux})
to a system of ordinary differential equations (ODEs) in
residual form. In the case where $p'=p$, we denote the residual
$\func{\rbm}{\Rbb^{N_\ubm}\times\Rbb^{N_\ubm}\times\Rbb^{N_\xbm}\times\Rbb^{N_\xbm}}{\Rbb^{N_\ubm}}$, 
which is defined as
\begin{equation}
 \rbm : (\mathring\wbm, \wbm, \mathring\ybm, \ybm) \mapsto
  \mbm \mathring\wbm + \fbm(\wbm, \ybm, \mathring\ybm),
\end{equation}
where $N_\ubm = \dim\Vcal_h^p$, $N_\xbm = \dim\Wcal_h$,
$\mbm\in\Rbb^{N_\ubm\times N_\ubm}$ is the mass matrix
associated with the test/trial space $\Vcal_h^p$, and
$\func{\fbm}{\Rbb^{N_\ubm}\times\Rbb^{N_\xbm}\times\Rbb^{N_\xbm}}{\Rbb^{N_\ubm}}$ 
is the algebraic form of the second and third terms in (\ref{eqn:claw-weak-numflux}).
In this notation, the standard DG discretization reads: given
$\func{\xbm}{[0,T]}{\Rbb^{N_\xbm}}$, find $\func{\ubm}{[0,T]}{\Rbb^{N_\ubm}}$
such that
\begin{equation} \label{eqn:claw-dg-semi}
 \rbm(\dot\ubm(t),\ubm(t),\dot\xbm(t),\xbm(t)) = \zerobold, \quad
 \ubm(0)=\bar\ubm, % \quad \xbm(0)=\bar\xbm,
\end{equation}
for all $t\in[0,T]$, where $\ubm$ is the time-dependent coefficients
of the DG solution, $\xbm$ is the time-dependent coefficients of
the domain mapping (nodal coordinates of the mesh), and
$\bar\ubm\in\Rbb^{N_\ubm}$ is the algebraic representation of
the initial condition $\bar{U}_X$; additionally, we define
$\bar\xbm\in\Rbb^{N_\xbm}$ as the initial condition for the nodal
coordinates, i.e., $\bar\xbm = \xbm(0)$.
Typically, the evolution of the mesh coordinates
$\xbm(t)$ is known analytically or governed by a dynamical
system (e.g., fluid-structure interaction); however, in this
work, it will be determined as the solution of an optimization
problem (after temporal discretization) such that discontinuities
are tracked over time.

Finally, we use the expansions in the nodal bases to define the
\emph{enriched residual}
$\func{\Rbm}{\Rbb^{N_\ubm}\times\Rbb^{N_\ubm}\times\Rbb^{N_\xbm}\times\Rbb^{N_\xbm}}{\Rbb^{N_\ubm'}}$
associated with a trial space of degree $p'$ as
\begin{equation}
 \Rbm : (\mathring\wbm, \wbm, \mathring\ybm, \ybm) \mapsto
  \Mbm \mathring\wbm + \Fbm(\wbm, \ybm, \mathring\ybm),
\end{equation}
where $N_\ubm' = \dim\Vcal_h^{p'}$,
$\Mbm\in\Rbb^{N_\ubm'\times N_\ubm}$ is the mass matrix
associated with the test space $\Vcal_h^{p'}$ and trial
space $\Vcal_h^p$, and
$\func{\Fbm}{\Rbb^{N_\ubm}\times\Rbb^{N_\xbm}\times\Rbb^{N_\xbm}}{\Rbb^{N_\ubm'}}$ 
is the algebraic form of the second and third terms in
(\ref{eqn:claw-weak-numflux}).
In this work, we take $p' = p +1$, but other choices are possible as well.
The enriched residual will be used in Section \ref{sec:optim}
to define the implicit tracking objective function.

\subsection{Numerical flux function}
\label{sec:disc:numflux}
The numerical flux function  corresponding to the reference flux 
($\numflux_X$, i.e. reference numerical flux) is a quantity that replaces the transformed flux 
dotted with the outward unit normal $(F_X \cdot N)$, 
as done from (\ref{eqn:claw-weak-elem}) to (\ref{eqn:claw-weak-elem-numflux}).
It is given by
\begin{equation} \label{eqn:claw-transf-numflux}
 \mapped\numflux(\mapped\stvc^+, \mapped\stvc^-, N; G, v) =
 \norm{\defjac\defgrad^{-T}N} \tilde\numflux(\stvc^+, \stvc^-, n; v),
\end{equation}
since by (\ref{eqn:transf-stvc})-(\ref{eqn:normal-transf}) we have
\begin{equation} \label{eqn:claw-transf-flux-dot-normal}
 \mapped\numflux \sim
 \mapped\flux \cdot N =
 \defjac(\flux - U \otimes v)\cdot\defgrad^{-T}N =
 \norm{\defjac\defgrad^{-T}N} (\flux - U \otimes v)\cdot n \sim
 \norm{\defjac\defgrad^{-T}N} \tilde\numflux,
\end{equation}
where arguments have been dropped for brevity. Here,
$\func{\tilde{\numflux}}{\Rbb^m\times\Rbb^m\times\Rbb^d\times\Rbb^d}{\Rbb^m}$
is a modified numerical flux which corresponds
to the modified flux function,
$\func{\tilde{F}}{\Rbb^m\times\Rbb^d}{\Rbb^{m\times d}}$,
which is related to the physical flux as
\begin{equation}
 \tilde{F}(U; v) = F(U) - U\otimes v.
\end{equation}
The modified flux accounts for the domain motion and
is obtained from the numerical flux function corresponding to 
the physical flux,
$\func{\numflux}{\Rbb^m\times\Rbb^m\times\Rbb^d}{\Rbb^m}$.
Notice that the Jacobian of the modified flux function only
differs from the Jacobian of the physical flux by a scale
multiple of the identity matrix (with scale factor $v\cdot n$),
which makes implementation of many numerical fluxes that
depend on the eigenvalue decomposition of the Jacobian
matrix, e.g., local Lax-Friedrichs, Roe, Vijayasundaram,
straightforward given the decomposition of the Jacobian
of the physical flux.

\begin{figure}
  \centering
  \input{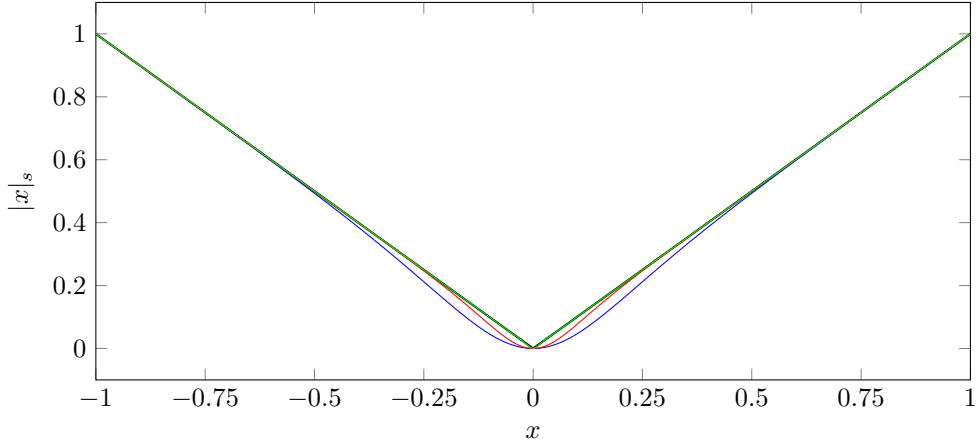}
  \caption{Smoothed absolute value function for $k = 5$ (\protect\bluesolidthinline), 
  	       $k = 10$ (\protect\redsolidthinline), and $k = 100$ (\protect\greensolidthinline).
  	       $k = \infty$ corresponds to the unsmoothed absolute value function (\protect\blacksolidthinline).}
  \label{fig:smooth-abs}
\end{figure}

For example, consider linear advection of a scalar field 
$\func{U}{\dom}{\Rbb}$ in a spatially varying direction
$\func{\beta}{\dom}{\Rbb^d}$ governed by a conservation law
of the form (\ref{eqn:claw-phys}) (see (\ref{eqn:advec}) in
Section~\ref{sec:num-exp:advec}) with physical flux function
\begin{equation} \label{eqn:advecflux}
 F_\adv(U) = U \beta^T
\end{equation}
with the corresponding upwind numerical flux
\begin{equation} \label{eqn:upwind}
 \Hcal_\up(U^+, U^-, n) =
 \begin{cases}
  (\beta\cdot n) U^+ & \text{ if } \beta\cdot n \geq 0 \\
  (\beta\cdot n) U^- & \text{ if } \beta\cdot n < 0.
 \end{cases}
\end{equation}
This can equivalently be written in terms of the absolute value function
$| \cdot |: \Rbb \rightarrow \Rbb_{\geq 0}$ as
\begin{equation} \label{eqn:upwindabs}
 \Hcal_\up(U^+, U^-, n) =
 0.5\left[(\beta\cdot n)(U^+ + U^-) + (U^+ - U^-)|\beta \cdot n | \right].
\end{equation}
Due to the linearity of the flux, the modified flux function is
\begin{equation}
 \tilde{F}_\text{adv}(U) = U(\beta-v)^T,
\end{equation}
which has an identical form as (\ref{eqn:advecflux}), albeit with a
modified velocity field. Thus, the corresponding modified upwind
numerical flux function is
\begin{equation}
 \tilde{\Hcal}_\text{up}(U^+,U^-,n; v) =
 0.5\left[((\beta - v)\cdot n)(U^+ + U^-) + (U^+ - U^-)|(\beta - v) \cdot n | \right].
\end{equation}
In \cite{zahr2020implicit}, we showed that it is advantageous, in terms of
solver performance, to have a numerical flux function that is smooth with
respect to variations in the normal $n$. Following the approach in
\cite{zahr2020implicit, zahr2020react}, we introduce a smoothed version
the upwind flux where the absolute value function is replaced with a
smoothed absolute value function 
\begin{equation} \label{eqn:upwindabs-sm}
 \Hcal_\up^s(U^+, U^-, n) =
 0.5\left[(\beta\cdot n)(U^+ + U^-) + (U^+ - U^-)|\beta \cdot n |_s \right],
\end{equation}
where $| \cdot |_s: \Rbb \rightarrow \Rbb_{\geq 0}$ is a 
smooth approximation to the absolute value given by
\begin{equation}
|\,\cdot\,|_s : x \mapsto x \tanh(kx)
\end{equation}
and $k$ is a smoothness parameter (Fig. \ref{fig:smooth-abs}). 
In this work, we use $k = 100$.
Previously, we had expressed numerical fluxes in terms of a
smoothed Heaviside function, but we find the smoothed 
absolute value function performs better numerically as we 
are no longer smoothing a discontinuous function.
Then, the smoothed version of the modified upwind
numerical flux is given by
\begin{equation} \label{eqn:upwindabs-sm-mod}
 \tilde{\Hcal}_\up^s(U^+, U^-, n; v) =
 0.5\left[((\beta - v)\cdot n)(U^+ + U^-) + (U^+ - U^-)|(\beta - v) \cdot n |_s \right].
\end{equation}

Similarly, for the Euler equations, we use Roe's flux \cite{roe1981approximate}
with the absolute value of the eigenvalues replaced by the corresponding
smoothed absolute value; a detailed derivation of the smoothed Roe flux
for the reacting Euler equations (single reaction) is included in Appendix
A of \cite{zahr2020react}.
%eqn:upwindabs is /camoapp/autog/mpl/eval_advec1v1d_nomap_core_inumflux_upwind.mpl
%eqn:upwindabs-sm is /camoapp/autog/mpl/eval_advec1v1d_nomap_core_inumflux_upwindsm.mpl 
%eqn:upwindabs-sm-mod is /camoapp/autog/mpl/eval_advec1v1d_alenogclc_core_inumflux_upwindsm.mpl

\subsection{High-order temporal discretization}
\label{sec:disc:high-order-time-disc}

Proceeding with the method of lines, we discretize the dynamical
system in (\ref{eqn:claw-dg-semi}) using a diagonally implicit
Runge-Kutta method (DIRK) to yield a sequence of algebraic
systems of equations. Unlike fully implicit Runge-Kutta methods,
an $s$-stage DIRK scheme has a Butcher tableau
$(A,b,c)\in\Rbb^{s\times s}\times\Rbb^s\times\Rbb^s$
where $A$ is lower triangular (Table~\ref{tab:dirk0}).
\begin{table}
 \caption{Butcher tableau for DIRK schemes}
 \label{tab:dirk0}
 \centering
 \renewcommand\arraystretch{1.2}
  \begin{tabular} { l | c c c }
 $\cbm$ & \Abm \\ \hline
    & $\bbm^T$
 \end{tabular}
 \qquad = \qquad
\renewcommand\arraystretch{1.2}
 \begin{tabular} { l | c c c }
 $c_1$ & $a_{11}$ & &  \\
 $\vdotswithin{=}$ & $\vdots$ & $\ddots$ & \\
 $c_s$ & $a_{s1}$ & $\dots$ & $a_{ss}$ \\ \hline
    & $b_1$ & $\dots$ & $b_s$
 \end{tabular} 
\end{table}
As a result, the $i$th stage only depends on the solution for
stages $1,\dots,i$, allowing the stages to be solved sequentially.
Because the nodal coordinates $\xbm$ are unknown in the proposed setting,
we first recast (\ref{eqn:claw-dg-semi}) as the following coupled system
of ODEs
\begin{equation} \label{eqn:claw-dg-semi1}
\begin{aligned}
 \mbm \dot\ubm(t) + \fbm(\ubm(t), \xbm(t), \bm\nu(t)) &= \zerobold \\
 \dot\xbm(t) - \bm\nu(t) &= \zerobold \\
 \ubm(0) &= \bar\ubm \\
 \xbm(0) &= \bar\xbm
\end{aligned}
\end{equation}
for $t\in[0,T]$,
where $\func{\bm{\nu}}{[0,T]}{\Rbb^{N_\xbm}}$ is the unknown nodal
velocity, defined as
\begin{equation}
 \bm{\nu} : t \mapsto \dot\xbm(t).
\end{equation}
Next, we partition the time interval $[0,T]$ into $N_T$ intervals
of equal size $\Delta t=T/N_T$ with endpoints $\{t_n\}_{n=0}^{N_T}$,
where $t_0=0$ and $t_n = t_{n-1} + \Delta t$ for $n=1,\dots,N_T$.
In this setting, an $s$-stage DIRK discretization of
(\ref{eqn:claw-dg-semi}) with Butcher tableau
$(A,b,c)$ reads:
for $n=1,\dots,N_T$ and $i=1,\dots,s$,
\begin{equation} \label{eqn:dirk0}
\begin{aligned}
 &\ubm_0 = \bar\ubm, \quad
 \ubm_{n+1} = \ubm_n + \sum_{j=1}^s b_i \kbm_{n,j}^\ubm, \quad
 \ubm_{n,i} = \ubm_n + \sum_{j=1}^i a_{ij} \kbm_{n,j}^\ubm \\
 &\xbm_0 = \bar\xbm, \quad
 \xbm_{n+1} = \xbm_n + \sum_{j=1}^s b_i \kbm_{n,j}^\xbm, \quad
 \xbm_{n,i} = \xbm_n + \sum_{j=1}^i a_{ij} \kbm_{n,j}^\xbm \\
 &\mbm \kbm_{n,i}^\ubm = -\Delta t\fbm\left(\ubm_{n,i},\xbm_{n,i},\bm\nu_{n,i}\right), \quad
 \kbm_{n,i}^\xbm = \Delta t\bm\nu_{n,i}, \quad \bm\nu_{n,i} = \bm\nu(t_n+c_i\Delta t)
\end{aligned}
\end{equation}
where $\ubm_0,\ubm_n,\ubm_{n,i},\kbm_{n,i}^\ubm\in\Rbb^{N_\ubm}$ and
$\xbm_0,\xbm_n,\xbm_{n,i},\kbm_{n,i}^\xbm,\bm\nu_{n,i}\in\Rbb^{N_\xbm}$
are implicitly defined as the solution of (\ref{eqn:dirk0});
$\ubm_n\approx\ubm(t_n)$ and $\xbm_n\approx\xbm(t_n)$
are the state and mesh approximation at each time step
$n=0,\dots,T$; $\ubm_{n,i}\approx\ubm(t_n+c_i\Delta t)$ and
$\xbm_{n,i}\approx\xbm(t_n+c_i\Delta t)$ are the state and
mesh approximations at each stage $i=1,\dots,s$ of each
time interval $n=1,\dots,T$; $\kbm_{n,i}^\ubm$ and $\kbm_{n,i}^\xbm$
are the solution and mesh stage updates,
and $\bm\nu_{n,i}=\bm\nu(t_n+c_i\Delta t)$ is the mesh velocity.
Because the mesh velocity $\bm\nu$ is unknown, the system in
(\ref{eqn:dirk0}) can neither be solved nor evaluated in residual form.

Following the work in \cite{froehle2015nonlinear}, we replace the
unknown velocity function $\bm\nu$ with a modified velocity function,
$\func{\tilde{\bm\nu}}{[0,T]}{\Rbb^{N_\xbm}}$, that ensures the
corresponding modified position
$\func{\tilde{\xbm}}{[0,T]}{\Rbb^{N_\xbm}}$, defined as the solution of
\begin{equation} \label{eqn:mod-mesheq}
 \dot{\tilde\xbm}(t) - \tilde{\bm\nu}(t) = \zerobold, \qquad
 \tilde\xbm(0) = \bar\xbm
\end{equation}
for $t\in[0,T]$, agree with the original stages $\xbm_{n,i}$
when discretized with the same DIRK method. Discretization
of (\ref{eqn:mod-mesheq}) using the DIRK scheme reads: for
$n=1,\dots,N_T$ and $i=1,\dots,s$,
\begin{equation} \label{eqn:dirk0p5}
 \tilde\xbm_0 = \bar\xbm, \quad
 \tilde\xbm_{n+1} = \tilde\xbm_n + \sum_{j=1}^s b_i \tilde\kbm_{n,j}^\xbm, \quad
 \tilde\xbm_{n,i} = \tilde\xbm_n + \sum_{j=1}^i a_{ij} \tilde\kbm_{n,j}^\xbm, \quad
 \tilde\kbm_{n,i}^\xbm = \Delta t\tilde{\bm\nu}_{n,i}, \quad \tilde{\bm\nu}_{n,i} = \tilde{\bm\nu}(t_n+c_i\Delta t),
\end{equation}
where $\tilde\xbm_0,\tilde\xbm_n,\tilde\xbm_{n,i},\tilde\kbm_{n,i}^\xbm,\tilde{\bm\nu}_{n,i}\in\Rbb^{N_\xbm}$ are implicitly defined as the solution of
(\ref{eqn:dirk0p5}); $\tilde\xbm_n\approx\tilde\xbm(t_n)$ is the modified mesh
approximation at each time step $n=0,\dots,T$;
$\tilde\xbm_{n,i}\approx\tilde\xbm(t_n+c_i\Delta t)$ is the modified mesh
approximation at each stage $i=1,\dots,s$ of each time interval $n=1,\dots,T$;
$\tilde\kbm_{n,i}^\xbm$ is the modified mesh stage update,
and $\tilde{\bm\nu}_{n,i}=\tilde{\bm\nu}(t_n+c_i\Delta t)$ is the modified
mesh velocity at each stage. Because the DIRK scheme only depends on the
value of $\tilde{\bm\nu}$ at $t_n+c_i\Delta t$, we only need to define 
its value at these points, i.e., $\tilde{\bm\nu}_{n,i}$; any smooth function
$\func{\zbm}{[0,T]}{\Rbb^{N_\xbm}}$ that satisfies
$\zbm(t_n+c_i\Delta t) = \tilde{\bm\nu}_{n,i}$ will lead to the
same DIRK solution and is a valid choice for $\tilde{\bm\nu}$. The
requirement that the modified and original mesh stages agree,
i.e., $\tilde\xbm_{n,i}=\xbm_{n,i}$ for $n=1,\dots,N_T$ and $i=1,\dots,s$,
leads to
\begin{equation}
 \xbm_{n,i} = \tilde\xbm_{n,i} =
 \tilde\xbm_n + \Delta t \sum_{j=1}^i a_{ij} \tilde{\bm\nu}_{n,j} =
 \xbm_n + \Delta t \sum_{j=1}^i a_{ij} \tilde{\bm\nu}_{n,j},
\end{equation}
where we used the definition of the stage $\tilde\xbm_{n,i}$ in
(\ref{eqn:dirk0p5}) and the relationship $\tilde\xbm_n = \xbm_n$
for $n=0,\dots,N_T$, which follows directly from the stage-consistent
requirement ($\tilde\xbm_{n,i}=\xbm_{n,i}$), the initial condition
($\tilde\xbm_0=\xbm_0=\bar\xbm$), and that the DIRK schemes used
to integrate (\ref{eqn:claw-dg-semi1}) and (\ref{eqn:mod-mesheq})
are the identical (same Butcher tableau). In the case where $A$ is full rank,
which is the case for fully and diagonally implicit Runge-Kutta methods,
this equation can be inverted to express the modified velocity stages
as a linear combination of the mesh stages
\begin{equation} \label{eqn:vel-stage-consist}
\tilde{\bm{\nu}}_{n,i}=\sum_{j=1}^i (\Abm^{-1})_{ij} \frac{\xbm_{n,j}-\xbm_n}{\Delta t}.
\end{equation}
This approach leads to a \emph{stage consistent} velocity approximation
\cite{froehle2015nonlinear} in the sense that the modified velocity
function is chosen such that it is consistent with the original mesh
positions at each stage.

%Following the work in \cite{froehle2015nonlinear}, we drop the
%requirement that $\bm\nu_{n,i}=\bm\nu(t_n+c_i\Delta t)$ and use
%(\ref{eqn:dirk0}) to relate the mesh and mesh velocity stages by
%\begin{equation}
% \xbm_{n,i} = \xbm_n + \Delta t \sum_{j=1}^i a_{ij} \bm\nu_{n,j}
%\end{equation}
%for $n=1,\dots,N_T$ and $i=1,\dots,s$. In the case where $A$ is full rank,
%which is the case for fully and diagonally implicit Runge-Kutta methods,
%this equation can be inverted to express the velocity stages as a linear
%combination of the mesh stages
%\begin{equation} \label{eqn:vel-stage-consist}
%\bm{\nu}_{n,i}=\sum_{j=1}^i (\Abm^{-1})_{ij} \frac{\xbm_{n,j}-\xbm_n}{\Delta t}.
%\end{equation}
%This approach leads to a \emph{stage consistent} velocity approximation
%\cite{froehle2015nonlinear}, i.e.,
%$\bm\nu_{n,i}=\tilde{\bm\nu}(t_n+c_i \Delta t)$, where
%$\func{\tilde{\bm\nu}}{[0,T]}{\Rbb^{N_\xbm}}$ is
%the modified velocity function that ensures the
%corresponding modified position
%$\func{\tilde{\xbm}}{[0,T]}{\Rbb^{N_\xbm}}$,
%defined as the solution of
%\begin{equation}
% \dot{\tilde\xbm}(t) - \tilde{\bm\nu}(t) = \zerobold, \qquad
% \tilde\xbm(0) = \bar\xbm
%\end{equation}
%for $t\in[0,T]$, agrees with stages $\xbm_{n,i}$
%when discretized with the same DIRK method.

With the stage consistent mesh velocity approximation, the
DIRK discretization in (\ref{eqn:dirk0}) becomes
\begin{equation} \label{eqn:dirk1}
\begin{aligned}
 &\ubm_0 = \bar\ubm, \quad
 \ubm_{n+1} = \ubm_n + \sum_{j=1}^s b_i \kbm_{n,j}^\ubm, \quad
 \ubm_{n,i} = \ubm_n + \sum_{j=1}^i a_{ij} \kbm_{n,j}^\ubm \\
 &\xbm_0 = \bar\xbm, \quad
 \xbm_{n+1} = \xbm_n + \sum_{j=1}^s b_i \kbm_{n,j}^\xbm, \quad
 \xbm_{n,i} = \xbm_n + \sum_{j=1}^i a_{ij} \kbm_{n,j}^\xbm \\
 &\mbm \kbm_{n,i}^\ubm = -\Delta t\fbm\left(\ubm_{n,i},\xbm_{n,i},\tilde{\bm\nu}_{n,i}\right), \quad
 \kbm_{n,i}^\xbm=\Delta t \tilde{\bm\nu}_{n,i}
\end{aligned}
\end{equation}
for $n=1,\dots,N_T$ and $i=1,\dots,s$, where the stage velocity
is defined in (\ref{eqn:vel-stage-consist}). Even though the unknown
velocity function has been eliminated using the stage consistent
velocity approximation, the new system (\ref{eqn:dirk1}) is
underdetermined because there are effectively $N_\ubm$ equations
in $N_\ubm+N_\xbm$ unknowns at a fixed stage. This will be resolved
by the optimization-based tracking formulation in the next section.

To close this section, we convert the modified DIRK system in
(\ref{eqn:dirk1}) to residual form
$\func{\rbm_{n,i}}{\Rbb^{N_\ubm}\times\Rbb^{N_\xbm}}{\Rbb^{N_\ubm}}$
at a fixed step $n\in\{1,\dots,N_T\}$ and stage $i\in\{1,\dots,s\}$ as
\begin{equation}
 \rbm_{n,i}: (\wbm,\ybm) \mapsto
   \mbm \xibold_{n,i}(\wbm) + \Delta t \fbm(\wbm, \ybm, \zetabold_{n,i}(\ybm)),
\end{equation}
where $\func{\xibold_{n,i}}{\Rbb^{N_\ubm}}{\Rbb^{N_\ubm}}$ maps the 
state stage ($\ubm_{n,i}$) to the corresponding stage update
($\kbm_{n,i}^\ubm$)
\begin{equation}
 \xibold_{n,i} : \wbm \mapsto
  (A^{-1})_{ii}(\wbm-\ubm_n) +
 \sum_{j=1}^{i-1} (A^{-1})_{ij}(\ubm_{n,j}-\ubm_n)
\end{equation}
and $\func{\zetabold_{n,i}}{\Rbb^{N_\xbm}}{\Rbb^{N_\xbm}}$ maps the
mesh stage ($\xbm_{n,i}$) to the corresponding stage-consistent velocity
($\tilde{\bm\nu}_{n,i}$) as
\begin{equation}
 \zetabold_{n,i} : \ybm \mapsto
  (A^{-1})_{ii}\frac{\ybm-\xbm_n}{\Delta t}+
 \sum_{j=1}^{i-1} (A^{-1})_{ij}\frac{\xbm_{n,j}-\xbm_n}{\Delta t}.
\end{equation}
Similarly, we define the corresponding fully discrete enriched residual
function $\func{\Rbm_{n,i}}{\Rbb^{N_\ubm}\times\Rbb^{N_\xbm}}{\Rbb^{N_\ubm'}}$
as
\begin{equation}
 \Rbm_{n,i} : (\wbm,\ybm) \mapsto
   \Mbm \xibold_{n,i}(\wbm) + \Delta t \Fbm(\wbm, \ybm, \zetabold_{n,i}(\ybm)),
\end{equation}
which we use in next section to define the objective function of
the implicit tracking optimization problem. Notice that for
$\Rbm_{n,i}$ only the \textit{spatial} test space is
enriched; the temporal discretization in $\rbm_{n,i}$
and $\Rbm_{n,i}$ is identical.

In this work, we consider three $L-$ stable DIRK
schemes with order of accuracy equal to the number of stages $s$;
Butcher tableaus given in Table~\ref{tab:DIRKbutcher}.
We denote the $k$-th order accurate DIRK scheme as DIRK$k$, i.e.
1st order DIRK is DIRK1, etc. We note that DIRK1 is equivalent to
the Backward Euler scheme.

\begin{table}
 \caption{Butcher tableau for DIRK1 (\textit{left}),
 DIRK2 (\textit{middle}), and DIRK3 (\textit{right}),
 where $\alpha = 1-\frac{1}{\sqrt{2}}$,
 $\beta=0.435866521508459$, $\gamma=-\frac{6\beta^2-16\beta+1}{4}$,
 and $\omega=\frac{6\beta^2-20\beta+5}{4}$.}
 \label{tab:DIRKbutcher}
\centering
\begin{minipage}{0.25\textwidth}
 \centering
 \begin{tabular} { r | c}
  1 & 1 \\ \hline
    & 1
 \end{tabular}
 \end{minipage}
 \begin{minipage}{0.25\textwidth}
 \centering
 \begin{tabular} { r | c c }
  $\alpha$ & $\alpha$ & 0 \\
   1    & 1-$\alpha$ & $\alpha$ \\ \hline
    & 1-$\alpha$ & $\alpha$
 \end{tabular}
 \end{minipage}
 \begin{minipage}{0.25\textwidth}
 \centering
 \begin{tabular} { r | c c c}
   $\beta$ & $\beta$ & & \\ 
   $\frac{1+ \beta}{2}$ & $\frac{1+ \beta}{2}-\beta$ & $\beta$ & \\
   $\gamma+\omega+\beta$ & $\gamma$ & $\omega$ & $\beta$ \\ \hline
   & $\gamma$ & $\omega$ & $\beta$
 \end{tabular}
 \end{minipage}
\end{table}

\section{Optimization formulation of $r$-adaptivity for implicit
         tracking of discontinuities}

\label{sec:optim}
In this section, we extend the high-order implicit shock
tracking framework from \cite{zahr2018optimization,zahr2020implicit} to
time-dependent problems using a method of lines approach (in contrast
to the space-time approach in \cite{corrigan2019moving, zahr2020implicit}).
Following the approach in \cite{zahr2020implicit}, we recast the fully discrete
conservation law as a PDE-constrained optimization problem over the 
discrete solution and mesh that aims to align element faces with discontinuities
at each Runge-Kutta stage for every timestep. We also review the associated
SQP solver for this optimization problem developed in \cite{zahr2020implicit}.

\subsection{Constrained optimization formulation}
\label{sec:optim:constr}

We formulate the problem of tracking discontinuities at
a given stage $i$ and timestep $n$ as a constrained optimization
problem over the PDE state and mesh stage that minimizes an
objective function, $\func{f_{n,i}}{\Rbb^{N_\ubm}\times\Rbb^{N_\xbm}}{\Rbb}$,
while enforcing the DG-DIRK discretization. That is, we define
$\ubm_{n,i}\in\Rbb^{N_\ubm}$ and $\xbm_{n, i}\in\Rbb^{N_\xbm}$ as
\begin{equation} \label{eqn:pde-opt-uns}
 (\ubm_{n,i},\xbm_{n,i}) \coloneqq
 \argmin_{\wbm\in\Rbb^{N_\ubm},\ybm\in\Rbb^{N_\xbm}} f_{n,i}(\wbm,\ybm)
 \quad\text{subject to:}\quad\rbm_{n,i}(\wbm,\ybm)=\zerobold.
\end{equation}
for $n=1,\dots,N_T$ and $i=1,\dots,s$. For a fixed time step $n$, 
once the stage states $\{\ubm_{n,i}\}_{i=1}^s$
and meshes $\{\xbm_{n,i}\}_{i=1}^s$ are computed, the state and mesh can
be advanced to the next time step ($\ubm_{n+1}$ and $\xbm_{n+1}$) using
the relationships in (\ref{eqn:dirk1}). The DIRK schemes considered in
this work (Table~\ref{tab:DIRKbutcher}) satisfy the property that
$A_{si} = b_i$ for $i=1,\dots,s$, which implies the state and mesh
at the final stage of time step $n$ are identical to the
state and mesh at time step $n+1$, i.e.,
$\ubm_{n+1} = \ubm_{n,s}$ and $\xbm_{n+1} = \xbm_{n,s}$.

The objective function is constructed such that
the solution of the PDE-constrained optimization problem is a feature-aligned
mesh by using the norm of the enriched DG-DIRK residual
\begin{equation} \label{eqn:obj-enrich-res-uns}
 f_{n, i}: (\wbm,\ybm) \mapsto
 \frac{1}{2}\norm{\Rbm_{n, i}(\wbm,\ybm)}_2^2.
\end{equation}
This follows on a large body of work that uses residual-based
error indicators to drive $h$-, $p$-, and $r$-adaptivity
\cite{fidkowski2011output}. This is the same objective function
used in our previous work for steady conservation laws, 
where it was shown to lead to a robust and reliable tracking framework;
for details, see \cite{zahr2020implicit}.

Unlike the steady case discussed in \cite{zahr2020implicit}, we do not
include the mesh quality term in the objective function for the
unsteady (method of lines) case. In the steady case, there is
no information about the shock location \textit{a priori}, which
usually requires significant deformation to the initial mesh to
align with shocks, necessitating the use of a mesh regularization term in the objective function.
However, in the context of timestepping, we have useful information
from the previous timestep to use as an initial guess for both the state
and mesh. We can combine this information along with the Rankine-Hugoniot
conditions  and a high-order mesh smoothing procedure to obtain an excellent
initial guess for the mesh at each Runge-Kutta stage
(Section \ref{sec:practical}).
In a sense, the additional time dimension allows the mesh regularization
to be decoupled from the implicit shock tracking procedure, which is one
advantage of the method of lines discretization over the space-time
formulation for solving unsteady problems.

\subsection{SQP solver for optimization-based discontinuity tracking}
\label{sec:optim:solver}

The Lagrangian of the optimization problem in (\ref{eqn:pde-opt-uns})
$\func{\Lcal}{\Rbb^{N_\ubm}\times\Rbb^{N_\xbm}\times\Rbb^{N_\ubm}}{\Rbb}$
takes the form
\begin{equation} \label{eqn:lagr}
 \Lcal(\ubm,\xbm,\lambdabold) = f(\ubm, \xbm)-\lambdabold^T\rbm(\ubm,\xbm),
\end{equation}
where $\lambdabold \in \Rbb^{N_\ubm}$ is a vector of Lagrange multipliers
associated with the DG constraint in (\ref{eqn:pde-opt-uns}).
We drop all the subscripts in this section for brevity. 
The first-order optimality, or Karush-Kuhn-Tucker (KKT), conditions state that the
$(\ubm^\star, \xbm^\star)$ is a first-order solution of the optimization
problem if there exists $\lambdabold^\star$ such that
\begin{equation} \label{eqn:kkt1}
\begin{gathered}
 \pder{f}{\ubm}(\ubm^\star,\xbm^\star)^T -
 \pder{\rbm}{\ubm}(\ubm^\star,\xbm^\star)^T\lambdabold^\star = \zerobold,
 \qquad
 \pder{f}{\xbm}(\ubm^\star,\xbm^\star)^T -
 \pder{\rbm}{\xbm}(\ubm^\star,\xbm^\star)^T\lambdabold^\star = \zerobold,
 \qquad
\rbm(\ubm^\star, \xbm^\star) = \zerobold.
\end{gathered}
\end{equation}
Since the DG Jacobian with respect to the state variables $\ubm$ is assumed
to be invertible, we define the estimate of the optimal Lagrange multiplier
$\func{\hat\lambdabold}{\Rbb^{N_\ubm}\times \Rbb^{N_\xbm}}{\Rbb^{N_\ubm}}$
such that the first equation ($\nabla_\ubm \Lcal = 0$) (adjoint equation) is
always satisfied
\begin{equation} \label{eqn:lagrmult}
 \hat\lambdabold(\ubm,\xbm) =
 \pder{\rbm}{\ubm}(\ubm,\xbm)^{-T}\pder{f}{\ubm}(\ubm,\xbm)^T.
\end{equation}
Then the optimality criteria becomes
\begin{equation}
 \cbm(\ubm^\star,\xbm^\star) \coloneqq
 \pder{f}{\xbm}(\ubm^\star,\xbm^\star)^T -
 \pder{\rbm}{\xbm}(\ubm^\star,\xbm^\star)^T\pder{\rbm^\star}{\ubm^\star}(\ubm^\star,\xbm^\star)^{-T}
 \pder{f}{\ubm}(\ubm^\star,\xbm^\star)^T = \zerobold, \qquad
 \rbm(\ubm^\star, \xbm^\star) = \zerobold
\end{equation}

Because the optimization problem (\ref{eqn:pde-opt-uns})
\textit{exactly matches the form} of the optimization problem for the
steady case \cite{zahr2020implicit}, we use the same SQP solver
that defines the sequences
$\{\ubm^{(k)}\}_{k=0}^\infty\subset\Rbb^{N_\ubm}$ and
$\{\xbm^{(k)}\}_{k=0}^\infty\subset\Rbb^{N_\xbm}$ as
\begin{equation} \label{eqn:sqpupdate}
 \ubm^{(k+1)} = \ubm^{(k)} + \alpha_{k+1} \Delta\ubm^{(k+1)}, \qquad
 \xbm^{(k+1)} = \xbm^{(k)} + \alpha_{k+1} \Delta\xbm^{(k+1)},
\end{equation}
for $k=0,1,\dots$, where $\alpha_{k+1}\in(0,1]$ is a step length
which can be determined by a line search procedure,
and $\Delta\ubm^{(k+1)}\in\Rbb^{N_\ubm}$ and
$\Delta\xbm^{(k+1)}\in\Rbb^{N_\xbm}$ are search directions.
At a given iteration $k$, the search directions
$\Delta\ubm^{(k)}$ and $\Delta\xbm^{(k)}$ are computed
simultaneously as the solution of a quadratic approximation
to the optimization problem in (\ref{eqn:pde-opt-uns}) with
a regularized Levenberg-Marquardt approximation of the Hessian;
for details, see \cite{zahr2020implicit}. A pair $(\ubm, \xbm)$
is considered a solution of (\ref{eqn:pde-opt-uns}) if
$\|\cbm(\ubm, \xbm)\| < \epsilon_1$ and
$\|\rbm(\ubm, \xbm)\| < \epsilon_2$ for tolerances
$\epsilon_1, \epsilon_2 > 0$.

\section{Practical considerations}
\label{sec:practical}

\subsection{Initial guess for optimization}
\label{sec:practical:MOLInitguess}

\begin{figure}
  \centering
 \begin{tikzpicture}

\begin{groupplot}[
    group style={ 
        group size=3 by 1,
        horizontal sep=0.5cm
    },
    enlargelimits=false,
    axis on top, axis equal image,
    xmin=0, xmax=1, ymin=0, ymax=1,
    ticks=none,
    width=0.45\textwidth,
]

\nextgroupplot
\addplot graphics [xmin=0, xmax=1, ymin=0, ymax=1] {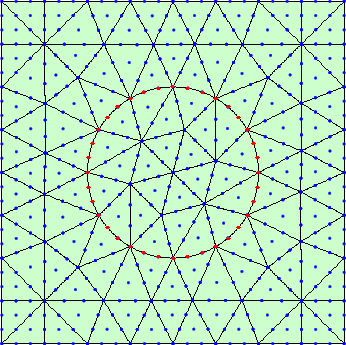};
%\draw[black, dashed] (axis cs:0.5,0.5) circle [radius=0.25];

\nextgroupplot
\addplot graphics [xmin=0, xmax=1, ymin=0, ymax=1] {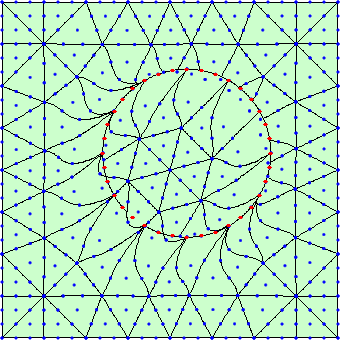};
%\draw[black, dashed] (axis cs:0.5,0.5) circle [radius=0.25];

\nextgroupplot
\addplot graphics [xmin=0, xmax=1, ymin=0, ymax=1] {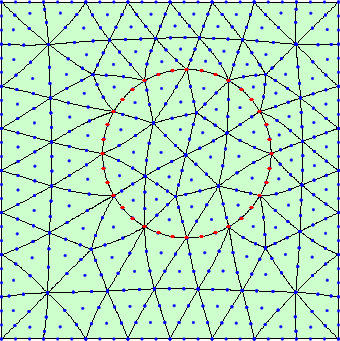};
%\draw[black, dashed] (axis cs:0.5,0.5) circle [radius=0.25];

\end{groupplot}
\end{tikzpicture}
 \caption{Shock nodes (\protect\reddottedline) (\emph{left}),
 	      advected by instantaneous shock speed (\emph{middle}) 
	      and smoothed mesh (\emph{right}).}
  \label{fig:mshsmooth}               
\end{figure}

The implicit shock tracking optimization problem (\ref{eqn:pde-opt-uns})
is non-convex and therefore the initial guess for the SQP solver
is critical to obtain a good solution. 
In the method of lines setting where implicit shock tracking is
performed at each time step, an obvious idea for the initial guess
is the converged mesh and solution from the previous time $t_n$.
In practice, this is too far off to attain good convergence
properties except for prohibitively small choices of timestep $\Delta t$. 
Instead, we employ an initial guess where we advect each node on the
discontinuity surface from the previous time by the instantaneous shock
speed determined by the Rankine-Hugoniot conditions  (Fig. \ref{fig:mshsmooth}). 
The remaining nodes are updated using a standard optimization-based 
mesh smoothing using the high-order mesh distortion metric 
$\func{r_K^\msh}{\Wcal_h}{\Rbb}$ developed in \cite{roca16distortion}, defined as
\begin{equation} \label{eqn:res-mshdist}
 r_\msh^K : \Qcal \mapsto
 \int_{K}
 \left(
  \frac{\norm{\nabla_X\Qcal}_F^2}{d(\det\nabla_X\Qcal)_+^{2/d}}
 \right)^2 \, dv.
\end{equation}
This provides a good initial guess for the shock-aligned mesh for 
each stage of the optimization problem. 
However, the advection of these shock nodes can result in a tangled mesh 
(Fig. \ref{fig:mshsmooth}, middle). We modify the mesh distortion metric 
(\ref{eqn:res-mshdist}) following the approach in \cite{gargallo2015optimization}
where the Jacobian $(\det\nabla_X\Qcal)$ is regularized to allow for the 
optimization procedure to recover from initially invalid configurations.

For the initial guess for the solution, we use the converged
physical solution ($U$) from the previous timestep.
However, since we applied our DG discretization to the transformed
reference conservation law (\ref{eqn:claw-ref}), we are solving for the
reference solution ($U_X$). Therefore, to use the physical solution
at time $t_n$ as the initial guess for time $t_{n+1}$, we multiply
it by the ratio of the Jacobian of the initial guess for the mapping
at time $t_{n+1}$ (advection based on Rankine-Hugoniot conditions and
smoothing) to the Jacobian of the converged mapping at time $t_n$.
More precisely, the initial guess for the $i$th 
Runge-Kutta stage in \textit{physical} space, denoted
$\check{U}_h$, is obtained by constant
extrapolation of the physical state at time $t$
\begin{equation}
\check{U}_h(x, t+c_i\Delta t) \coloneqq U_h(x, t) = g_h(X, t)^{-1} U_{X_h}(X, t).
\end{equation}
Let $\tilde{\Gcal} \in \Wcal_h$ be the mapping to the configuration
obtained by the Rankine-Hugoniot-based procedure and associated smoothing,
then the initial guess for the $i$th Runge-Kutta stage in the
\textit{reference} domain, denoted $\check{U}_{X_h}$, is
\begin{equation}
\check{U}_{X_h}(X, t+c_i\Delta t) \coloneqq \tilde{g}(X, t + c_i\Delta t) \check{U}_h(x,t+c_i\Delta t) = \frac{\tilde{g}(X, t + c_i\Delta t)}{g_h(X, t)}U_{X_h}(X, t).
\end{equation}
Future work to enhance the  robustness of these initial guesses might
consider more advanced methods to approximate the shock velocity
\cite{rawat2010high} and higher order extrapolation-based estimates
for the solution.

\subsection{Line search}
\label{sec:practical:lsrch}
In the steady case, the initial mesh is generated independently of the
\textit{a priori} unknown shock location, and as such 
can begin quite far from alignment with the shock. 
Element collapses were required to remove small and poorly shaped elements 
that arose between iterations as the coarse high-order mesh deformed to align with the shock.
Because of this, we chose the step length parameter $\alpha_{k+1}$ of (\ref{eqn:sqpupdate})
based off a backtracking line search procedure based on the $\ell_1$ merit function
in our previous work \cite{zahr2020implicit}. 

In this work, we simply take $\alpha_k = 1$ and do not require any topology changes to the mesh.
This is possible because we are able to construct high-quality initial guesses 
for the optimization problem in the unsteady case that already contain information 
about the shock (unlike the steady case).
In Section \ref{sec:num-exp}, we demonstrate this choice leads to
{\color{blue}{rapid convergence}} of the SQP solver.
%Note this does not quite correspond to Newton's method, due to the Gauss-Newton
%approximation to the Hessian used (Section 4.4-4.6 of \cite{zahr2020implicit}), but
%in practice behaves very similarly, which we demonstrate in some of our
%numerical experiments in Section \ref{sec:num-exp}.

\section{Numerical experiments}
\label{sec:num-exp}
In this section, we introduce three inviscid unsteady conservation laws
and demonstrate the tracking framework on five problems with
discontinuous solutions of varying difficulty. We also 
demonstrate high-order convergence in time and present 
the convergence behavior of the SQP solver. 

\subsection{Linear advection}
\label{sec:num-exp:advec}
The first problem we consider is linear advection of a scalar
quantity $U$ through a domain $\Omega \subset \Rbb^d$
\begin{equation} \label{eqn:advec}
\begin{aligned}
 \pder{}{t}U(x,t) + \pder{}{x_j}\left(U(x, t)\beta_j(x)\right) & = 0 \quad &&\text{for}~x\in\Omega,~t\in[0,T] \\
 U(x, 0) &= \bar{U}(x) &&\text{for}~x\in\Omega,
\end{aligned}
\end{equation}
where $\func{U}{\dom\times[0,T]}{\Rbb}$ is the conserved quantity
implicitly defined as the solution of (\ref{eqn:advec}),
$\func{\beta}{\Omega}{\Rbb^d}$ is the flow direction, and 
$\func{\bar{U}}{\Omega}{\Rbb}$ is the initial condition.

\subsubsection{1D spatially varying advection}
\label{sec:num-exp:advec:1dspadvec}
\begin{figure}
 \centering
 \input{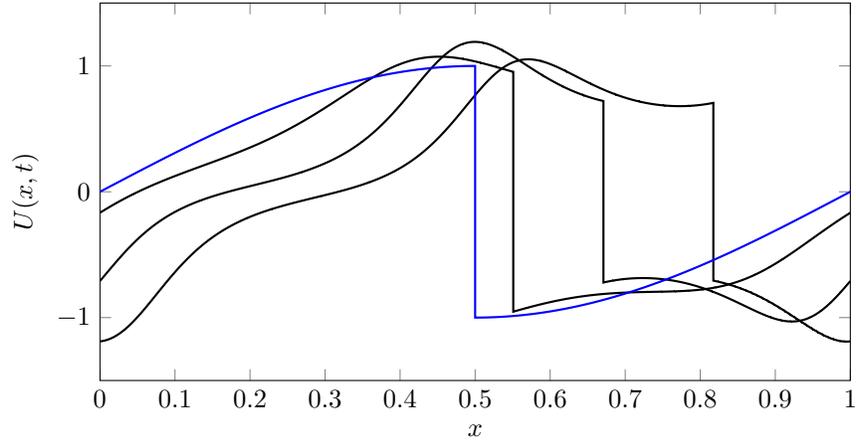}
 \caption{Method of lines solution of the one-dimensional, spatially varying advection equation with $p = 4, q = 1$.
 	      Initial condition $\bar{U}(x)$ (\ref{line:advec1d_mol_IC}) and tracking solution 
	      at times $t = 0.05, 0.15, 0.25$ (\ref{line:advec1d_mol_slice}).}
 \label{fig:advec1d_mol_slice} 
\end{figure}

\begin{figure}
 \centering
 %%%This was adapted from advec2d_beta0_c2s2_nel8x4_porder0x2_conv.tikz and the corresponding python script from the 2020 JCP paper

\begin{tikzpicture}
\begin{groupplot} [
group style={group size = 2 by 2, horizontal sep = 2.3cm, vertical sep = 2.1cm}]
\nextgroupplot[width=0.4\textwidth, ymin=1e-12, ytick={1e-8, 1e-10, 1e-12}, xlabel={Timestep ($n$)}, ymax=2e-8, xmax=25, mark repeat={5}, xmin=0, ymode=log, ylabel={$\norm{\rbm(\ubm_n, \xbm_n)}$}]
\addplot [solid, thick, color=black, mark options={solid, thin}, mark=*, mark size=1.5, color=black]
coordinates {
(1 ,0.000000000335174396586396)
(2 ,0.000000002087545266306)
(3 ,0.00000000536082073096822)
(4 ,3.85855552147893E-12)
(5 ,1.6760193674263E-12)
(6 ,2.7244644644291E-12)
(7 ,4.18903100669669E-12)
(8 ,6.16068444537169E-12)
(9 ,8.69842056044088E-12)
(10 ,1.17834934802951E-11)
(11 ,1.52556671455074E-11)
(12 ,1.87486290977658E-11)
(13 ,2.16702466047925E-11)
(14 ,2.32797024466527E-11)
(15 ,2.28979762486441E-11)
(16 ,2.02071335840075E-11)
(17 ,1.55038168765805E-11)
(18 ,9.7636649073965E-12)
(19 ,4.43562164238621E-12)
(20 ,0.00000000473037342893424)
(21 ,0.000000000616946237769889)
(22 ,1.0697972285281E-11)
(23 ,7.14394505647698E-12)
(24 ,1.33194246060934E-11)
(25 ,1.9979516606387E-11)
};\label{line:advec_conv:values}

\addplot [dashed, thick, color=pink, mark=none, forget plot]
coordinates {
( 0,  1e-08)
( 110,  1e-08)
};;\label{line:advec_conv:tolerances}

\nextgroupplot[width=0.4\textwidth, ymin=1e-10, ytick={1e-6, 1e-8, 1e-10}, xlabel={Timestep ($n$)}, ymax=2e-6, xmax=25, mark repeat={5}, xmin=0, ymode=log, ylabel={$\norm{\cbm(\ubm_n, \xbm_n)}$}]
\addplot [solid, thick, color=black, mark options={solid, thin}, mark=*, mark size=1.5, color=black, forget plot]
coordinates {
(1 ,0.000000141670033328629)
(2 ,0.000000355674976795883)
(3 ,0.000000570700276414204)
(4 ,0.00000000119497023960212)
(5 ,0.00000000139405748074058)
(6 ,0.00000000162716678163783)
(7 ,0.00000000197408333297913)
(8 ,0.00000000252977974987561)
(9 ,0.000000003260981417007)
(10 ,0.00000000418276208680514)
(11 ,0.00000000527613138050266)
(12 ,0.00000000666285658174505)
(13 ,0.00000000839268371918135)
(14 ,0.0000000104609394308277)
(15 ,0.0000000127319633730919)
(16 ,0.000000014954930800246)
(17 ,0.0000000168473185238204)
(18 ,0.0000000180770293249921)
(19 ,0.0000000182790786531757)
(20 ,0.000000340166854496211)
(21 ,0.000000122606173463013)
(22 ,0.0000000108296726096796)
(23 ,0.00000000666539112282702)
(24 ,0.00000000293154861056597)
(25 ,0.00000000058491039602027)
};

\addplot [dashed, thick, color=pink, mark=none, forget plot]
coordinates {
( 0,  1e-06)
( 25,  1e-06)
};;\label{line:advec_conv:tolerances}

\nextgroupplot[width=0.4\textwidth, ymin=1e-7, ytick={1e-5, 1e-6, 1e-7}, xlabel={Timestep ($n$)}, ymax=1e-5, xmax=25, mark repeat={5}, xmin=0, ymode=log, ylabel={$\norm{\Rbm(\ubm_n, \xbm_n)}$}]
\addplot [solid, thick, color=black, mark options={solid, thin}, mark=*, mark size=1.5, color=black, forget plot]
coordinates {
(1 ,0.000000570255847544434)
(2 ,0.000000447950400765092)
(3 ,0.000000531232820065309)
(4 ,0.000000554324202225737)
(5 ,0.000000620669470104049)
(6 ,0.000000667005071430721)
(7 ,0.0000007618265278685)
(8 ,0.00000087289030349285)
(9 ,0.00000101143486012022)
(10 ,0.00000117069175479227)
(11 ,0.00000134079956198848)
(12 ,0.00000154278475460106)
(13 ,0.00000177044709252831)
(14 ,0.00000202231623649255)
(15 ,0.00000228143723029493)
(16 ,0.00000252498014267703)
(17 ,0.00000273303962884285)
(18 ,0.00000288037856083301)
(19 ,0.00000293650405255051)
(20 ,0.00000287116315705891)
(21 ,0.00000265323313128724)
(22 ,0.00000224014634442032)
(23 ,0.00000168020603326436)
(24 ,0.00000102239821241509)
(25 ,0.000000511425643285399)
};

\nextgroupplot[
width=0.4\textwidth,
xlabel={Timestep ($n$)},
ymax=3,
xmax=25,
ylabel={No. of SQP steps},
xmin=0,
ymin=1,
mark repeat={5}]
\addplot [solid, thick,mark options={solid, thin}, mark=*, mark size=1.5, color=blue]
coordinates {
(1, 2)
(2, 2)
(3, 2)
(4, 2)
(5, 2)
(6, 2)
(7, 2)
(8, 2)
(9, 2)
(10, 2)
(11, 2)
(12, 2)
(13, 2)
(14, 2)
(15, 2)
(16, 2)
(17, 2)
(18, 2)
(19, 2)
(20, 2)
(21, 2)
(22, 2)
(23, 2)
(24, 2)
(25, 2)
};\label{line:advec_newton:stage1}

\addplot [dashed, thick, mark options={solid, thin}, mark=square, mark size=1.5, color=red]
coordinates {
(1, 2)
(2, 2)
(3, 2)
(4, 2)
(5, 2)
(6, 2)
(7, 2)
(8, 2)
(9, 2)
(10, 2)
(11, 2)
(12, 2)
(13, 2)
(14, 2)
(15, 2)
(16, 2)
(17, 2)
(18, 2)
(19, 2)
(20, 2)
(21, 2)
(22, 2)
(23, 2)
(24, 2)
(25, 2)
};\label{line:advec_newton:stage2}

\addplot [dotted, thick, mark options={solid, thin}, mark=diamond, mark size=1.5, color=magenta]
coordinates {
(1, 2)
(2, 2)
(3, 2)
(4, 3)
(5, 3)
(6, 3)
(7, 3)
(8, 3)
(9, 3)
(10, 3)
(11, 3)
(12, 3)
(13, 3)
(14, 3)
(15, 3)
(16, 3)
(17, 3)
(18, 3)
(19, 3)
(20, 3)
(21, 2)
(22, 3)
(23, 3)
(24, 3)
(25, 3)
};\label{line:advec_newton:stage3}

\end{groupplot}\end{tikzpicture}
 \caption{Final converged value of the constraint (\textit{top left}),
          optimality condition (\textit{top right}), and objective
          function (\textit{bottom left}) for the solution to 1D spatially varying advection at 
          each timestep (\ref{line:advec_conv:values})
          and the specified tolerances (\ref{line:advec_conv:tolerances}),
          and the number of SQP steps needed at each stage of each time
          step for convergence (\textit{bottom right});
          stage 1: (\ref{line:advec_newton:stage1}),
	        stage 2 (\ref{line:advec_newton:stage2}),
	        stage 3 (\ref{line:advec_newton:stage3}).}
 \label{fig:advec_sqp}
\end{figure}
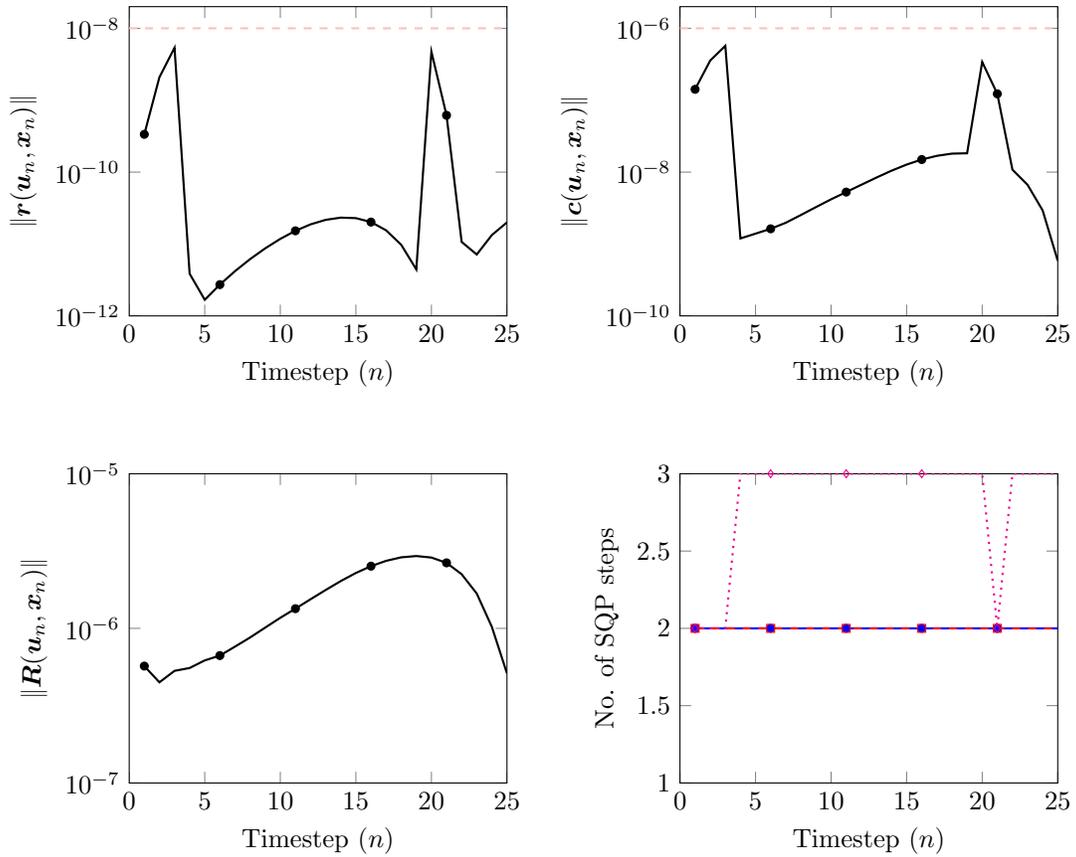

As a simple benchmark problem to demonstrate the capabilities
and performance of the unsteady shock tracking framework, we
consider the advection equation in one spatial dimension
with a spatially varying advection field $\beta: \Omega \rightarrow \Rbb$
\begin{equation}
 \beta : x \mapsto 1 + \frac{1}{2}\sin^2(2\pi x)
\end{equation}
with initial condition $\bar{U}: \Omega \rightarrow \Rbb$:
\begin{equation}
\bar{U}(x) = 
\begin{cases} 
      \sin(\pi x) &  x \leq 0.5 \\
      \sin(\pi(x-1)) & x > 0.5
   \end{cases}
\end{equation}
and periodic boundary conditions.
We initialize with a mesh of the reference domain $\Omega_0=[0,1]$,
which we construct such that an element interface lies at the
initial shock location ($x=0.5$), i.e., the shock in the initial
condition is tracked. The shock tracking solution is computed
using a DG discretization on this equispaced mesh with $20$ elements of
degree $p=4$, $q=1$ and a DIRK$3$ temporal discretization with
$25$ time steps with final time $T = 0.25$ (Figure~\ref{fig:advec1d_mol_slice}).
The SQP solver is used with tolerances $\epsilon_1 = 10^{-6}$
and $\epsilon_2 = 10^{-8}$. At each timestep $n+1$, we are able
obtain a solution to within the specified tolerances, 
which corresponds to the solution of the optimization problem
(\ref{eqn:pde-opt-uns}) at time step $n$ and stage $s$,
and only need 2-3 iterations of the SQP solver to do so (Figure~\ref{fig:advec_sqp}).
Note that the $s$th stage tends to be more difficult for the SQP solver --
one possible explanation is because each stage of the DIRK method
has low stage order (first order), making it easier on the SQP solver
compared to  the $s$th stage, which is in fact the high-order solution
because of our  specific choice of DIRK schemes.
We acknowledge this is an easy example for the SQP solver, and
will demonstrate similar behavior for a more challenging example
in the following sections. 

\begin{figure}
 \centering
 \input{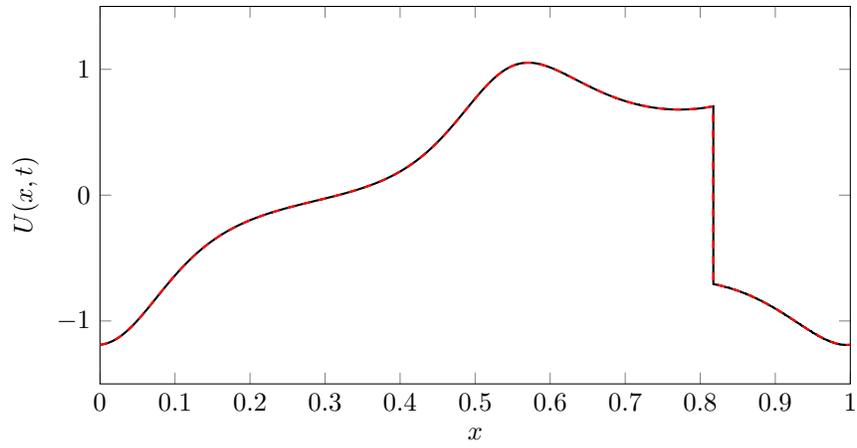}
 \caption{Comparison of shock tracking (\ref{line:advec1d_mol_slice}) and reference (\ref{line:advec1d_mol_ref}) solutions at $T = 0.25$.}
 \label{fig:advec1d_refcomp}
\end{figure}

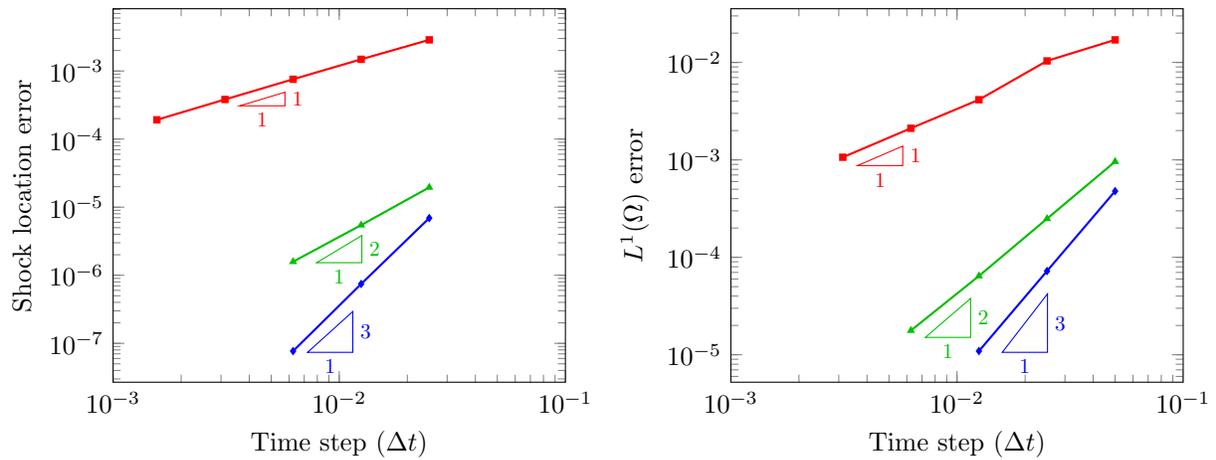
\begin{figure}
 \centering
 \begin{tikzpicture}

\begin{groupplot}[
    group style={ 
        group size=2 by 1,
        horizontal sep=2.2cm
    },
]

\nextgroupplot[width=0.46\textwidth, xmin=0.001, xmax=0.1, xlabel={Time step ($\Delta t$)}, xmode = log, ymode=log, ylabel={Shock location error}]
\addplot [red, solid, thick, mark=square*, mark size=1, mark options={solid}]  
coordinates {
(0.025, 0.002867396583290)
(0.0125, 0.001484845533757)
(0.00625, 0.000756218466129)
(0.003125, 0.000381793494416)
(0.0015625,  0.000191914526769)
};\label{line:DIRK1}
\addplot [green!75!black, solid, thick, mark=triangle*, mark size=1, mark options={solid}] 
coordinates {
(0.025, 0.000019512146133)
(0.0125, 0.000005455831985)
(0.00625, 0.000001581123337)
};\label{line:DIRK2}
\addplot [blue, solid, thick, mark=diamond*, mark size=1, mark options={solid}] 
coordinates {
(0.025, 0.000006912729842)
(0.0125, 0.000000744830044)
(0.00625, 0.000000077002040)
};\label{line:DIRK3}

\logLogSlopeTriangle{0.38}{0.1}{0.74}{1}{red};
\logLogSlopeTriangle{0.55}{0.1}{0.32}{2}{green!75!black};
\logLogSlopeTriangle{0.53}{0.1}{0.08}{3}{blue};

\nextgroupplot[width=0.46\textwidth, xmin=0.001, xmax=0.1, xlabel={Time step ($\Delta t$)}, xmode = log, ymode=log, ylabel={$L^1(\Omega)$ error}]
\addplot [red, solid, thick, mark=square*, mark size=1, mark options={solid}] 
coordinates {
(0.05, 0.017056278173930)
(0.025, 0.010373642096652)
(0.0125, 0.0041434086369029)
(0.00625, 0.002111140512689)
(0.003126, 0.001064293155348)
};\label{line:DIRK1}
\addplot [green!75!black, solid, thick, mark=triangle*, mark size=1, mark options={solid}] 
coordinates {
(0.05, 0.000960559082526)
(0.025,0.000249230808895)
(0.0125, 0.000064366724884)
(0.00625, 0.000017722463395)
};\label{line:DIRK2}
\addplot [blue, solid, thick, mark=diamond*, mark size=1, mark options={solid}] 
coordinates {
(0.05, 0.000477580093943)
(0.025, 0.000072370197858)
(0.0125, 0.000010914821411)
};\label{line:DIRK3}

\logLogSlopeTriangle{0.38}{0.1}{0.58}{1}{red};
\logLogSlopeTriangle{0.53}{0.1}{0.12}{2}{green!75!black};
\logLogSlopeTriangle{0.70}{0.1}{0.08}{3}{blue};

\end{groupplot}
\end{tikzpicture}
 \caption{Temporal convergence of the DIRK1 (\ref{line:DIRK1}), DIRK2 (\ref{line:DIRK2}), and DIRK3 (\ref{line:DIRK3})
 	      schemes for the $L^1$ error of the shock location (\textit{left}) and the solution (\textit{right}) 
	      for the spatially varying advection equation.}
 \label{fig:advec1d_highorder_conv}
\end{figure}

We obtain a reference solution of (\ref{eqn:advec}) at $T = 0.25$ using
the method of characteristics; the characteristic equations and corresponding
solution are integrated using classical RK4 with 10,000 timesteps.
The shock tracking solution compares well to the reference
solution at the final time $T = 0.25$ (Figure~\ref{fig:advec1d_refcomp}).
We use this reference solution to demonstrate high-order convergence 
in time of both the $L^1$ error of the solution at the final time and
the shock location,  verifying the design order of accuracy of the
DIRK$k$ schemes even in the presence of a discontinuity
(Figure~\ref{fig:advec1d_highorder_conv}).

%\begin{figure}
% \centering
% \input{py/advec_conv.tikz}
% \caption{Convergence history of the constraint (\textit{left}),
%          optimality condition (\textit{middle}), and objective
%          function (\textit{right}) of 1D spatially varying advection at 
%          each timestep (\ref{line:advec_conv:upwind})
%          and the specified tolerances (\ref{line:advec_conv:tolerances}).}          
% \label{fig:advec_conv}
%\end{figure}
%
%\begin{figure}
% \centering
% \input{py/advec_newton.tikz}
% \caption{Number of SQP steps needed at each timestep for 
%	      stage 1 (\ref{line:advec_newton:stage1}),
%	      stage 2 (\ref{line:advec_newton:stage2}),
%	      and stage 3 (\ref{line:advec_newton:stage3})
%	      for the spatially varying advection equation.
%        {\color{red}{I think this would be better on the same plot.}}}
% \label{fig:advec_newton} 
%\end{figure}
%
%\begin{figure}
% \centering
% \input{py/advec_newton_sameplt.tikz}
% \caption{Number of SQP steps needed at each timestep for 
%	      stage 1 (\ref{line:advec_newton:stage1}),
%	      stage 2 (\ref{line:advec_newton:stage2}),
%	      and stage 3 (\ref{line:advec_newton:stage3})
%	      for the spatially varying advection equation.
%        {\color{red}{I think this would be better on the same plot.}}}
% \label{fig:advec_newton2} 
%\end{figure}

\subsection{Time-dependent, inviscid Burgers' equation}
\label{sec:num-exp:burg}
The time-dependent, inviscid Burgers' equation
governs nonlinear advection of a scalar quantity through the
domain $\Omega \subset \Rbb^d$
\begin{equation} \label{eqn:burg}
\begin{aligned}
 \pder{}{t}U(x,t) + \pder{}{x_j}\left(\frac{1}{2}U(x,t)^2\beta_j\right) &= 0 \quad &&\text{for}~x\in\Omega,~t\in[0,T] \\
 U(x, t) &= 0 &&\text{for}~x\in\partial\Omega,~t\in[0,T] \\
 U(x, 0) &= \bar{U}(x) &&\text{for}~x\in\Omega,
\end{aligned}
\end{equation}
where $\func{U}{\dom\times[0,T]}{\Rbb}$ is the conserved quantity
implicitly defined as the solution of (\ref{eqn:burg}),
$\beta\in\Rbb^d$ is the flow direction, and 
$\func{\bar{U}}{\Omega}{\Rbb}$ is the initial condition.
We investigate this problem in both one and two dimensions.
\subsubsection{1D Burgers' equation}
\label{sec:num-exp:burg1d}

\begin{figure}
 \centering
 \input{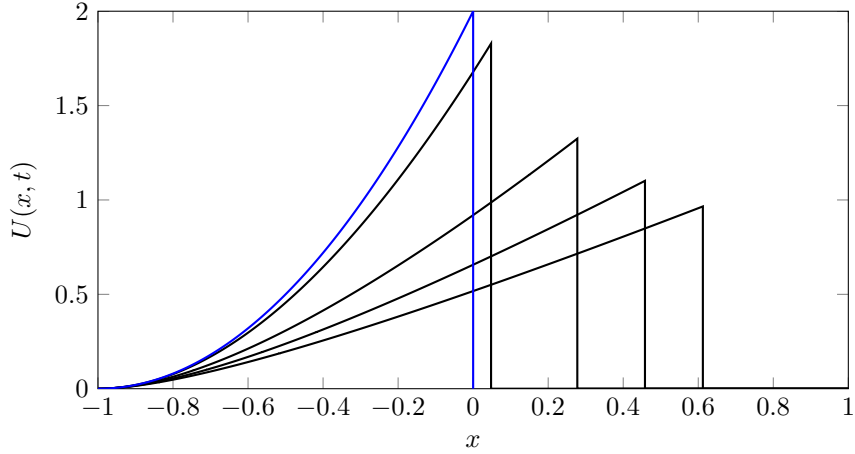}
 \caption{Method of lines solution of the one-dimensional, inviscid
 Burgers' equation with $p = 4, q = 1$, including initial condition
 $\bar{U}(x)$ (\ref{line:burg1d_mol_IC}) and tracking solution 
 at times $t = 0.05, 0.35, 0.65, 0.95$ (\ref{line:burg1d_mol_slice}).}
 \label{fig:burg1d_mol_slice}
\end{figure}

For Burgers' equation in one spatial dimension,
we take $\beta=1$ and consider the initial condition $\func{\bar{U}}{\Omega}{\Rbb}$
\begin{equation}
 \bar{U} : x \mapsto 2(x+1)^2(1-H(x)),
\end{equation}
where $\func{H}{\Rbb}{\{0,1\}}$ is the Heaviside function.
Our approach only requires a mesh of the reference domain $\Omega_0=[-1,1]$,
which we construct such that an element interface lies at the
initial shock location ($x=0$), i.e., the shock in the initial
condition is tracked. The shock tracking solution is computed
using a DG discretization on this equispaced mesh with $20$ elements of
degree $p=4$, $q=1$ and a DIRK$3$ temporal discretization with
$20$ time steps with final time $T = 1$ (Figure~\ref{fig:burg1d_mol_slice}).
%The method of lines solution
%(Figure \ref{fig:burg1dMOLsoln}) agrees well with the corresponding
%slice of the $p = q = 3$ solution obtained by the space-time formulation.
%\begin{figure}
% \centering
%  \includegraphics[width=0.45\textwidth]{img/burg1dMOLfig1.png}
%  \includegraphics[width=0.45\textwidth]{img/burg1dMOLfig2.png}
% \caption{Method of lines solution of one-dimensional, inviscid Burgers'
% equation with $p = 3$, $q = 1$.  Initial condition $\bar{U}(x)$ (\textit{left})
% and solution at $T = 1$ (\textit{right}).}
% \label{fig:burg1dMOLsoln}
%\end{figure}

\subsubsection{2D Burgers' equation}
\label{sec:num-exp:burg2d}

\begin{figure}
 \centering
  \includegraphics[width=0.45\textwidth]{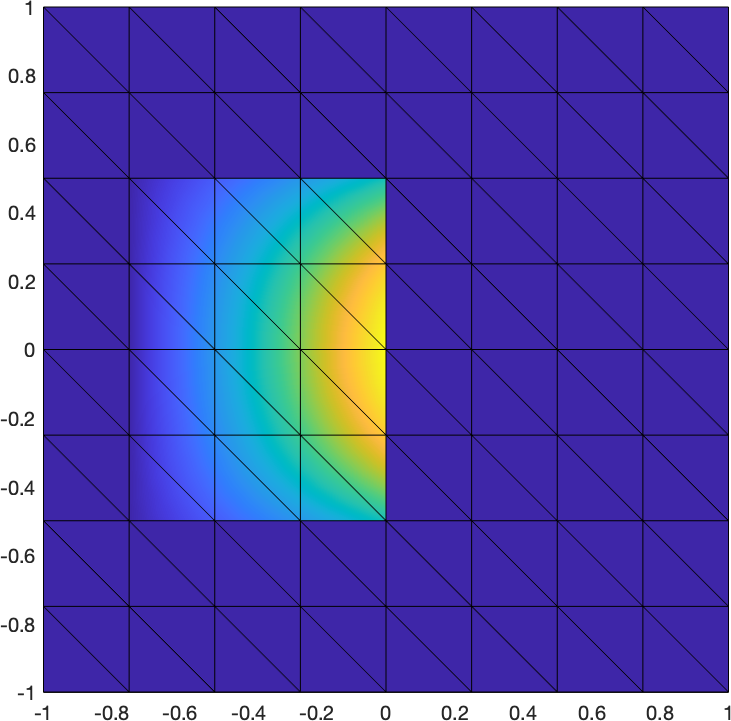}
  \includegraphics[width=0.45\textwidth]{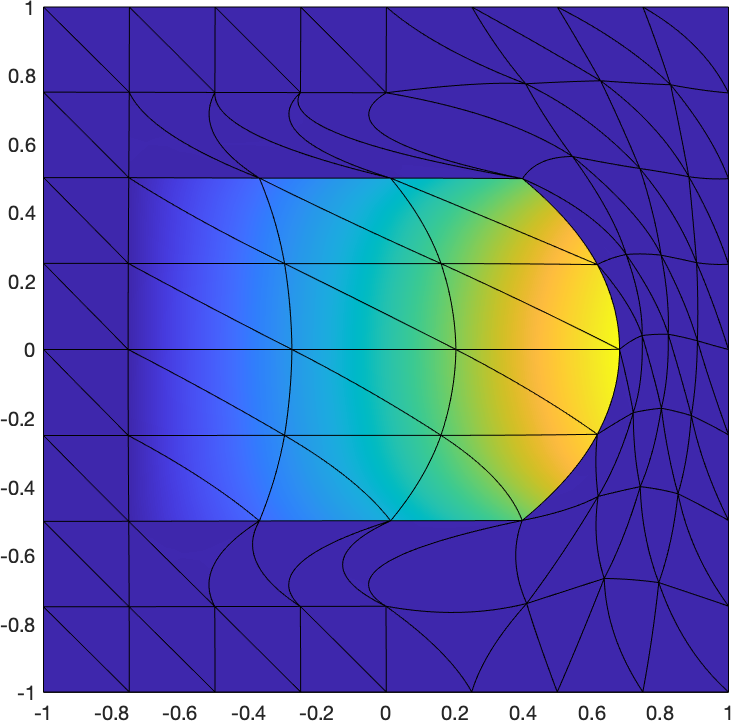}
 \caption{Method of lines solution of two-dimensional, inviscid Burgers' equation with $p = q = 2$. 
 Initial condition $\bar{U}(x)$ (\textit{left}) and solution at $T = 2$ (\textit{right})}
 \label{fig:burg2dMOLsoln}
 \begin{tikzpicture}
\begin{axis}[
   hide axis, scale only axis,
   height=0pt, width=0pt,
   colormap={parula}{rgb255=(62,38,168) rgb255=(62,39,172) rgb255=(63,40,175) rgb255=(63,41,178) rgb255=(64,42,180) rgb255=(64,43,183) rgb255=(65,44,186) rgb255=(65,45,189) rgb255=(66,46,191) rgb255=(66,47,194) rgb255=(67,48,197) rgb255=(67,49,200) rgb255=(67,50,202) rgb255=(68,51,205) rgb255=(68,52,208) rgb255=(69,53,210) rgb255=(69,55,213) rgb255=(69,56,215) rgb255=(70,57,217) rgb255=(70,58,220) rgb255=(70,59,222) rgb255=(70,61,224) rgb255=(71,62,225) rgb255=(71,63,227) rgb255=(71,65,229) rgb255=(71,66,230) rgb255=(71,68,232) rgb255=(71,69,233) rgb255=(71,70,235) rgb255=(72,72,236) rgb255=(72,73,237) rgb255=(72,75,238) rgb255=(72,76,240) rgb255=(72,78,241) rgb255=(72,79,242) rgb255=(72,80,243) rgb255=(72,82,244) rgb255=(72,83,245) rgb255=(72,84,246) rgb255=(71,86,247) rgb255=(71,87,247) rgb255=(71,89,248) rgb255=(71,90,249) rgb255=(71,91,250) rgb255=(71,93,250) rgb255=(70,94,251) rgb255=(70,96,251) rgb255=(70,97,252) rgb255=(69,98,252) rgb255=(69,100,253) rgb255=(68,101,253) rgb255=(67,103,253) rgb255=(67,104,254) rgb255=(66,106,254) rgb255=(65,107,254) rgb255=(64,109,254) rgb255=(63,110,255) rgb255=(62,112,255) rgb255=(60,113,255) rgb255=(59,115,255) rgb255=(57,116,255) rgb255=(56,118,254) rgb255=(54,119,254) rgb255=(53,121,253) rgb255=(51,122,253) rgb255=(50,124,252) rgb255=(49,125,252) rgb255=(48,127,251) rgb255=(47,128,250) rgb255=(47,130,250) rgb255=(46,131,249) rgb255=(46,132,248) rgb255=(46,134,248) rgb255=(46,135,247) rgb255=(45,136,246) rgb255=(45,138,245) rgb255=(45,139,244) rgb255=(45,140,243) rgb255=(45,142,242) rgb255=(44,143,241) rgb255=(44,144,240) rgb255=(43,145,239) rgb255=(42,147,238) rgb255=(41,148,237) rgb255=(40,149,236) rgb255=(39,151,235) rgb255=(39,152,234) rgb255=(38,153,233) rgb255=(38,154,232) rgb255=(37,155,232) rgb255=(37,156,231) rgb255=(36,158,230) rgb255=(36,159,229) rgb255=(35,160,229) rgb255=(35,161,228) rgb255=(34,162,228) rgb255=(33,163,227) rgb255=(32,165,227) rgb255=(31,166,226) rgb255=(30,167,225) rgb255=(29,168,225) rgb255=(29,169,224) rgb255=(28,170,223) rgb255=(27,171,222) rgb255=(26,172,221) rgb255=(25,173,220) rgb255=(23,174,218) rgb255=(22,175,217) rgb255=(20,176,216) rgb255=(18,177,214) rgb255=(16,178,213) rgb255=(14,179,212) rgb255=(11,179,210) rgb255=(8,180,209) rgb255=(6,181,207) rgb255=(4,182,206) rgb255=(2,183,204) rgb255=(1,183,202) rgb255=(0,184,201) rgb255=(0,185,199) rgb255=(0,186,198) rgb255=(1,186,196) rgb255=(2,187,194) rgb255=(4,187,193) rgb255=(6,188,191) rgb255=(9,189,189) rgb255=(13,189,188) rgb255=(16,190,186) rgb255=(20,190,184) rgb255=(23,191,182) rgb255=(26,192,181) rgb255=(29,192,179) rgb255=(32,193,177) rgb255=(35,193,175) rgb255=(37,194,174) rgb255=(39,194,172) rgb255=(41,195,170) rgb255=(43,195,168) rgb255=(44,196,166) rgb255=(46,196,165) rgb255=(47,197,163) rgb255=(49,197,161) rgb255=(50,198,159) rgb255=(51,199,157) rgb255=(53,199,155) rgb255=(54,200,153) rgb255=(56,200,150) rgb255=(57,201,148) rgb255=(59,201,146) rgb255=(61,202,144) rgb255=(64,202,141) rgb255=(66,202,139) rgb255=(69,203,137) rgb255=(72,203,134) rgb255=(75,203,132) rgb255=(78,204,129) rgb255=(81,204,127) rgb255=(84,204,124) rgb255=(87,204,122) rgb255=(90,204,119) rgb255=(94,205,116) rgb255=(97,205,114) rgb255=(100,205,111) rgb255=(103,205,108) rgb255=(107,205,105) rgb255=(110,205,102) rgb255=(114,205,100) rgb255=(118,204,97) rgb255=(121,204,94) rgb255=(125,204,91) rgb255=(129,204,89) rgb255=(132,204,86) rgb255=(136,203,83) rgb255=(139,203,81) rgb255=(143,203,78) rgb255=(147,202,75) rgb255=(150,202,72) rgb255=(154,201,70) rgb255=(157,201,67) rgb255=(161,200,64) rgb255=(164,200,62) rgb255=(167,199,59) rgb255=(171,199,57) rgb255=(174,198,55) rgb255=(178,198,53) rgb255=(181,197,51) rgb255=(184,196,49) rgb255=(187,196,47) rgb255=(190,195,45) rgb255=(194,195,44) rgb255=(197,194,42) rgb255=(200,193,41) rgb255=(203,193,40) rgb255=(206,192,39) rgb255=(208,191,39) rgb255=(211,191,39) rgb255=(214,190,39) rgb255=(217,190,40) rgb255=(219,189,40) rgb255=(222,188,41) rgb255=(225,188,42) rgb255=(227,188,43) rgb255=(230,187,45) rgb255=(232,187,46) rgb255=(234,186,48) rgb255=(236,186,50) rgb255=(239,186,53) rgb255=(241,186,55) rgb255=(243,186,57) rgb255=(245,186,59) rgb255=(247,186,61) rgb255=(249,186,62) rgb255=(251,187,62) rgb255=(252,188,62) rgb255=(254,189,61) rgb255=(254,190,60) rgb255=(254,192,59) rgb255=(254,193,58) rgb255=(254,194,57) rgb255=(254,196,56) rgb255=(254,197,55) rgb255=(254,199,53) rgb255=(254,200,52) rgb255=(254,202,51) rgb255=(253,203,50) rgb255=(253,205,49) rgb255=(253,206,49) rgb255=(252,208,48) rgb255=(251,210,47) rgb255=(251,211,46) rgb255=(250,213,46) rgb255=(249,214,45) rgb255=(249,216,44) rgb255=(248,217,43) rgb255=(247,219,42) rgb255=(247,221,42) rgb255=(246,222,41) rgb255=(246,224,40) rgb255=(245,225,40) rgb255=(245,227,39) rgb255=(245,229,38) rgb255=(245,230,38) rgb255=(245,232,37) rgb255=(245,233,36) rgb255=(245,235,35) rgb255=(245,236,34) rgb255=(245,238,33) rgb255=(246,239,32) rgb255=(246,241,31) rgb255=(246,242,30) rgb255=(247,244,28) rgb255=(247,245,27) rgb255=(248,247,26) rgb255=(248,248,24) rgb255=(249,249,22) rgb255=(249,251,21) },
   colorbar horizontal,
   point meta min=0.000000e+00, point meta max=1.000000e+00,
   colorbar style={width=10cm, xtick={0.000000e+00,2.500000e-01,5.000000e-01,7.500000e-01,1.000000e+00}}
]
\addplot [draw=none] coordinates {(0,0)};
\end{axis}
\end{tikzpicture}  
\end{figure}
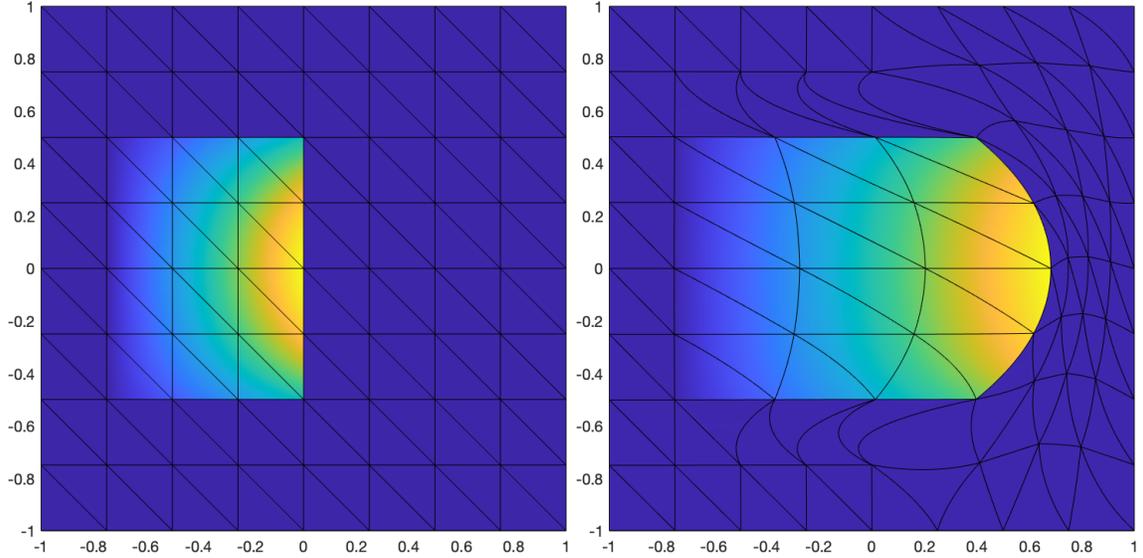

We now consider Burgers' equation in two spatial dimensions.
We take $\Omega_0=[-1,1]^2$ as our two-dimensional spatial domain,
$\beta=(1,0)$ as the flow direction, and 
$\func{\bar{U}}{\Omega}{\Rbb}$ as the initial condition, defined as
\begin{equation} \label{eqn:burg2d-mol-initcond}
 \bar{U} : (x_1, x_2) \mapsto
\begin{cases} 
      (0.5-2(x_2^2-0.25))
      \left(\frac{4}{3}(x_1 + 0.75)\right) &
      x\in\Omega_\square \\
      0 & \text{elsewhere}, 
   \end{cases}
\end{equation}
where $\Omega_\square\coloneqq [-0.75, 0]\times[-0.5,0.5]$.
The problem is constructed such that the initially straight
shock curves over time, which is tracked by the high-order
mesh. The shock tracking solution is computed using a DG discretization
on a mesh with $128$ simplex elements of degree $p=2$, $q=2$ and a DIRK$3$
temporal discretization with $40$ time steps with final time $T=2$
(Figure~\ref{fig:burg2dMOLsoln}).
The mesh smoothing procedure described in
Section~\ref{sec:practical:MOLInitguess} is important here to
maintain high-quality elements as the shock moves across the domain.

\subsection{Unsteady, compressible Euler equations}
\label{sec:num-exp:euler}

The Euler equations govern the flow of an inviscid, compressible fluid
through a domain $\Omega \subset \Rbb^d$
\begin{equation} \label{eqn:euler}
\begin{split}
\pder{}{t}\rho(x,t) + \pder{}{x_j}\left(\rho(x,t) v_j(x,t)\right) &= 0 \\
\pder{}{t}\left(\rho(x,t)v_i(x,t)\right) +
\pder{}{x_j}\left(\rho(x,t) v_i(x,t)v_j(x,t)+P(x,t)\delta_{ij}\right) &= 0 \\
\pder{}{t}\left(\rho(x,t)E(x,t)\right) +
\pder{}{x_j}\left(\left[\rho(x,t)E(x,t)+P(x,t)\right]v_j(x,t)\right) &= 0
\end{split}
\end{equation}
for all $x\in\Omega$, $t\in[0,T]$, $i=1,\dots,d$ and summation is implied
over the repeated index $j=1,\dots,d$, where
$\func{\rho}{\Omega\times (0, T)}{\Rbb_+}$ is
the density of the fluid, $\func{v_i}{\Omega\times (0, T)}{\Rbb}$ for
$i = 1, \dots, d$ is the velocity of the fluid in the $x_i$ direction,
and $\func{E}{\Omega\times (0, T)}{\Rbb_{+}}$ is the total energy of
the fluid, implicitly defined as the solution of (\ref{eqn:euler}).
For a calorically ideal fluid, the pressure of the fluid,
$\func{P}{\Omega\times (0, T)}{\Rbb_{+}}$, is related to the energy
via the ideal gas law
\begin{equation}
 P = (\gamma-1)\left(\rho E - \frac{\rho v_i v_i}{2}\right),
\end{equation}
where $\gamma\in\Rbb_{+}$ is the ratio of specific heats.
By combining the density, momentum, and energy into a vector of
conservative variables $\func{U}{\Omega\times[0,T]}{\Rbb^{d+2}}$,
defined as
\begin{equation}
 U : (x, t) \mapsto
 \begin{bmatrix}
  \rho(x,t) \\ \rho(x,t)v(x,t) \\ \rho(x,t)E(x,t)
 \end{bmatrix}
\end{equation}
the Euler equations are a conservation law of the form (\ref{eqn:claw-phys}).
We investigate the shock tracking framework on two benchmark examples:
the Shu-Osher problem and a blast wave problem in 2D.

\subsubsection{Shu-Osher Problem}
\label{sec:num-exp:euler:shuosher}

\begin{figure}
 \centering
 \includegraphics[width=0.75\textwidth]{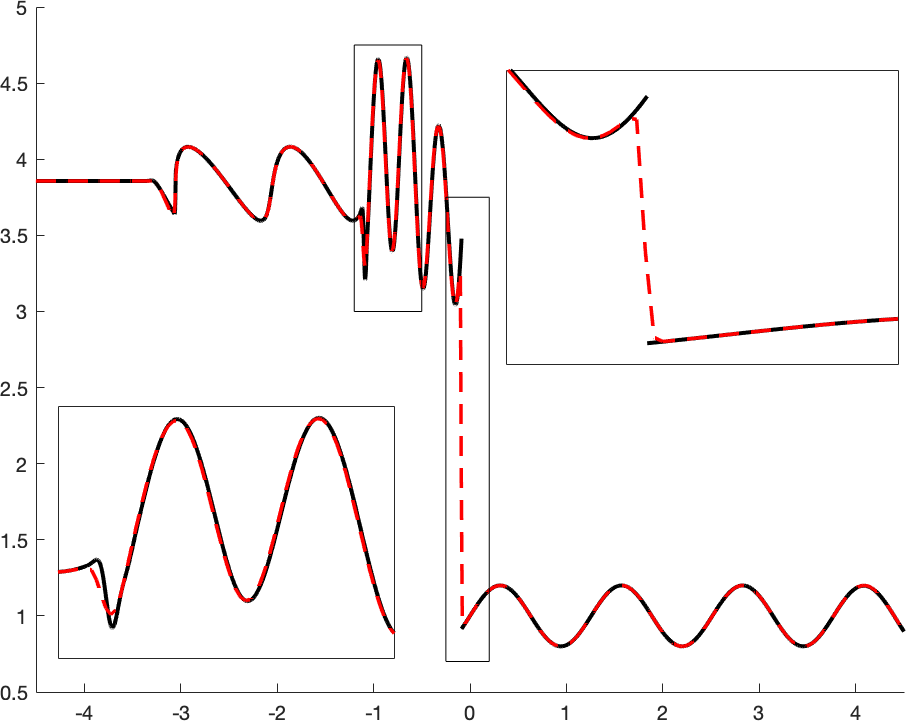}
 \caption{Density at $T = 1.1$ of Shu-Osher problem for the
 	        reference (\protect\reddashedline) and shock
          tracking (\protect\blacksolidline) solutions.} 
 \label{fig:shuosher}
\end{figure}

The Shu-Osher problem \cite{shu1989efficient} is a one-dimensional
idealization of shock-turbulence  interaction where a Mach 3 shock
moves into a field with a small sinusoidal density disturbance. The
flow domain is $\Omega_0 = [-4.5, 4.5]$, and the initial condition 
is given in terms of the density, velocity, and pressure as
\begin{equation}
\begin{aligned}
 \rho(x) &= \begin{cases} 3.857143 & x<-4 \\ 1+0.2\sin(5x) & x \geq -4 \end{cases} \\
 v(x) &= \begin{cases} 2.629369 & x<-4 \\ 0 & x \geq -4 \end{cases} \\
 P(x) &= \begin{cases} 10.3333 & x<-4 \\ 1 & x \geq -4, \end{cases}
\end{aligned}
\end{equation}
and the density, velocity, and pressure are prescribed at $x=-4.5$ and the
velocity is prescribed at $x=4.5$ (values can be read from the initial
condition). The shock tracking solution is computed using a DG
discretization on a mesh with $288$ elements of degree $p=4$, $q=1$,
half of which are equispaced to the left of the initial shock location from 
$[-4.5, -4]$, and the other half equispaced from $[-4, 4.5]$.
The temporal discretization is done by the DIRK$3$ 
method with $110$ time steps with final time $T = 1.1$.
The final time is chosen such that waves trailing behind the
primary shock do not steepen into shock waves; shock formation will be the
subject of future work. In Figure~\ref{fig:shuosher}, we
present the shock tracking solution along with a reference solution
computed using a fifth-order WENO method with 200 elements
and temporal integration via RK4 with $110$ timesteps
\cite{shu1989efficient}. The shock tracking solution actually overshoots the reference
solution at the formation of the trailing waves, which suggests the
reference solution is being overly dissipated by the WENO scheme
(\emph{left inset}). The shock is perfectly represented by the
aligned mesh in the shock tracking solution compared to the
reference (\emph{right inset}).

The SQP solver is used with tolerances $\epsilon_1 = 10^{-4}$
and $\epsilon_2 = 10^{-8}$, which is easily attained at each Runge-Kutta
stage and timestep and in very few SQP iterations (Figure~\ref{fig:shuosher_sqp}).
We note that the objective function is steadily growing because we are
using a fixed amount of resolution to represent an increasingly 
oscillatory numerical solution. This is not an issue of any other fundamental
concern, and can be addressed with adaptive element collapses and refinement,
which is straightforward in the one-dimensional case.

The Shu-Osher problem is ideal for implicit shock tracking because it starts
off with a well-defined shock that maintains its topology over time 
and the flow away from the shock is smooth. 
At present, one could imagine combining this shock tracking approach
with a shock capturing method, where strong well-defined shocks
could be tracked and weaker shocks which form over time could be captured \cite{johnsen2010assessment}.
Future work will focus on the simulation of the full Shu-Osher problem
to additionally track the smaller shocks which form as the trailing waves steepen.

\begin{figure}
 \centering
 \input{py/shuosher_sqp.tikz}
 \caption{Final converged value of the constraint (\textit{top left}),
          optimality condition (\textit{top right}), and objective
          function (\textit{bottom left}) for the solution to the Shu-Osher problem at
          each timestep (\ref{line:shuosher_conv:values})
          and the specified tolerances (\ref{line:shuosher_conv:tolerances}),
          and the number of SQP steps needed at each stage of each time
          step for convergence (\textit{bottom right});
          stage 1: (\ref{line:shuosher_newton:stage1}),
          stage 2 (\ref{line:shuosher_newton:stage2}),
          stage 3 (\ref{line:shuosher_newton:stage3}).}
 \label{fig:shuosher_sqp}
\end{figure}

\subsubsection{Blast Wave}
\label{sec:num-exp:euler:blast}

In the final example, we consider a spherical blast wave problem
featuring a strong, radially expanding shock. 
Our problem is inspired by the classic Sedov problem, 
which is an idealized model of the self-similar evolution of a cylindrical 
(or spherical) blast wave starting from a large total energy deposited
at a single point placed into an otherwise homogenous medium of a uniform
ambient density with negligible pressure. 
%Practically speaking numerical people initialize for a small finite radius and spread
%out the pressure due to the difficulty with representing a single point mass
%and pick a pressure negligibly small so the expansion velocity
%of the blast wave is much higher than the speed of sound in the ambient medium.
The Sedov problem is a workhorse verification test for traditional shock capturing
methods since there is an analytical solution \cite{sedov1993similarity} available 
in one-, two-, and three-dimensions for comparison. It is primarily used 
to test geometrical concerns, such as the ability to track a curved shock
and maintain spherical symmetry.

\begin{figure}
 \centering
  \includegraphics[width=0.45\textwidth]{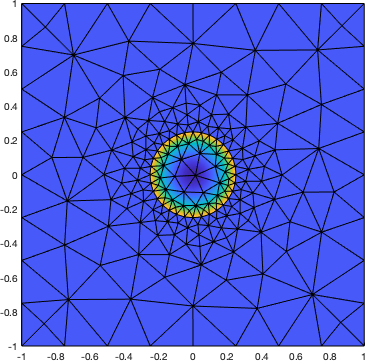}
  \includegraphics[width=0.45\textwidth]{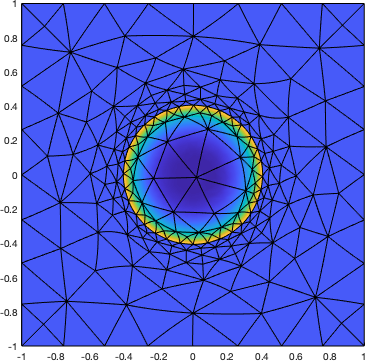}
 \caption{Method of lines solution of two-dimensional blast wave problem with $p = q = 2$. 
 Initial condition $\bar{U}(x)$ (\textit{left}) and solution at $T = 0.14$ (\textit{right})}
 \label{fig:blastwaveMOLsoln}
 \begin{tikzpicture}
\begin{axis}[
   hide axis, scale only axis,
   height=0pt, width=0pt,
   colormap={parula}{rgb255=(62,38,168) rgb255=(62,39,172) rgb255=(63,40,175) rgb255=(63,41,178) rgb255=(64,42,180) rgb255=(64,43,183) rgb255=(65,44,186) rgb255=(65,45,189) rgb255=(66,46,191) rgb255=(66,47,194) rgb255=(67,48,197) rgb255=(67,49,200) rgb255=(67,50,202) rgb255=(68,51,205) rgb255=(68,52,208) rgb255=(69,53,210) rgb255=(69,55,213) rgb255=(69,56,215) rgb255=(70,57,217) rgb255=(70,58,220) rgb255=(70,59,222) rgb255=(70,61,224) rgb255=(71,62,225) rgb255=(71,63,227) rgb255=(71,65,229) rgb255=(71,66,230) rgb255=(71,68,232) rgb255=(71,69,233) rgb255=(71,70,235) rgb255=(72,72,236) rgb255=(72,73,237) rgb255=(72,75,238) rgb255=(72,76,240) rgb255=(72,78,241) rgb255=(72,79,242) rgb255=(72,80,243) rgb255=(72,82,244) rgb255=(72,83,245) rgb255=(72,84,246) rgb255=(71,86,247) rgb255=(71,87,247) rgb255=(71,89,248) rgb255=(71,90,249) rgb255=(71,91,250) rgb255=(71,93,250) rgb255=(70,94,251) rgb255=(70,96,251) rgb255=(70,97,252) rgb255=(69,98,252) rgb255=(69,100,253) rgb255=(68,101,253) rgb255=(67,103,253) rgb255=(67,104,254) rgb255=(66,106,254) rgb255=(65,107,254) rgb255=(64,109,254) rgb255=(63,110,255) rgb255=(62,112,255) rgb255=(60,113,255) rgb255=(59,115,255) rgb255=(57,116,255) rgb255=(56,118,254) rgb255=(54,119,254) rgb255=(53,121,253) rgb255=(51,122,253) rgb255=(50,124,252) rgb255=(49,125,252) rgb255=(48,127,251) rgb255=(47,128,250) rgb255=(47,130,250) rgb255=(46,131,249) rgb255=(46,132,248) rgb255=(46,134,248) rgb255=(46,135,247) rgb255=(45,136,246) rgb255=(45,138,245) rgb255=(45,139,244) rgb255=(45,140,243) rgb255=(45,142,242) rgb255=(44,143,241) rgb255=(44,144,240) rgb255=(43,145,239) rgb255=(42,147,238) rgb255=(41,148,237) rgb255=(40,149,236) rgb255=(39,151,235) rgb255=(39,152,234) rgb255=(38,153,233) rgb255=(38,154,232) rgb255=(37,155,232) rgb255=(37,156,231) rgb255=(36,158,230) rgb255=(36,159,229) rgb255=(35,160,229) rgb255=(35,161,228) rgb255=(34,162,228) rgb255=(33,163,227) rgb255=(32,165,227) rgb255=(31,166,226) rgb255=(30,167,225) rgb255=(29,168,225) rgb255=(29,169,224) rgb255=(28,170,223) rgb255=(27,171,222) rgb255=(26,172,221) rgb255=(25,173,220) rgb255=(23,174,218) rgb255=(22,175,217) rgb255=(20,176,216) rgb255=(18,177,214) rgb255=(16,178,213) rgb255=(14,179,212) rgb255=(11,179,210) rgb255=(8,180,209) rgb255=(6,181,207) rgb255=(4,182,206) rgb255=(2,183,204) rgb255=(1,183,202) rgb255=(0,184,201) rgb255=(0,185,199) rgb255=(0,186,198) rgb255=(1,186,196) rgb255=(2,187,194) rgb255=(4,187,193) rgb255=(6,188,191) rgb255=(9,189,189) rgb255=(13,189,188) rgb255=(16,190,186) rgb255=(20,190,184) rgb255=(23,191,182) rgb255=(26,192,181) rgb255=(29,192,179) rgb255=(32,193,177) rgb255=(35,193,175) rgb255=(37,194,174) rgb255=(39,194,172) rgb255=(41,195,170) rgb255=(43,195,168) rgb255=(44,196,166) rgb255=(46,196,165) rgb255=(47,197,163) rgb255=(49,197,161) rgb255=(50,198,159) rgb255=(51,199,157) rgb255=(53,199,155) rgb255=(54,200,153) rgb255=(56,200,150) rgb255=(57,201,148) rgb255=(59,201,146) rgb255=(61,202,144) rgb255=(64,202,141) rgb255=(66,202,139) rgb255=(69,203,137) rgb255=(72,203,134) rgb255=(75,203,132) rgb255=(78,204,129) rgb255=(81,204,127) rgb255=(84,204,124) rgb255=(87,204,122) rgb255=(90,204,119) rgb255=(94,205,116) rgb255=(97,205,114) rgb255=(100,205,111) rgb255=(103,205,108) rgb255=(107,205,105) rgb255=(110,205,102) rgb255=(114,205,100) rgb255=(118,204,97) rgb255=(121,204,94) rgb255=(125,204,91) rgb255=(129,204,89) rgb255=(132,204,86) rgb255=(136,203,83) rgb255=(139,203,81) rgb255=(143,203,78) rgb255=(147,202,75) rgb255=(150,202,72) rgb255=(154,201,70) rgb255=(157,201,67) rgb255=(161,200,64) rgb255=(164,200,62) rgb255=(167,199,59) rgb255=(171,199,57) rgb255=(174,198,55) rgb255=(178,198,53) rgb255=(181,197,51) rgb255=(184,196,49) rgb255=(187,196,47) rgb255=(190,195,45) rgb255=(194,195,44) rgb255=(197,194,42) rgb255=(200,193,41) rgb255=(203,193,40) rgb255=(206,192,39) rgb255=(208,191,39) rgb255=(211,191,39) rgb255=(214,190,39) rgb255=(217,190,40) rgb255=(219,189,40) rgb255=(222,188,41) rgb255=(225,188,42) rgb255=(227,188,43) rgb255=(230,187,45) rgb255=(232,187,46) rgb255=(234,186,48) rgb255=(236,186,50) rgb255=(239,186,53) rgb255=(241,186,55) rgb255=(243,186,57) rgb255=(245,186,59) rgb255=(247,186,61) rgb255=(249,186,62) rgb255=(251,187,62) rgb255=(252,188,62) rgb255=(254,189,61) rgb255=(254,190,60) rgb255=(254,192,59) rgb255=(254,193,58) rgb255=(254,194,57) rgb255=(254,196,56) rgb255=(254,197,55) rgb255=(254,199,53) rgb255=(254,200,52) rgb255=(254,202,51) rgb255=(253,203,50) rgb255=(253,205,49) rgb255=(253,206,49) rgb255=(252,208,48) rgb255=(251,210,47) rgb255=(251,211,46) rgb255=(250,213,46) rgb255=(249,214,45) rgb255=(249,216,44) rgb255=(248,217,43) rgb255=(247,219,42) rgb255=(247,221,42) rgb255=(246,222,41) rgb255=(246,224,40) rgb255=(245,225,40) rgb255=(245,227,39) rgb255=(245,229,38) rgb255=(245,230,38) rgb255=(245,232,37) rgb255=(245,233,36) rgb255=(245,235,35) rgb255=(245,236,34) rgb255=(245,238,33) rgb255=(246,239,32) rgb255=(246,241,31) rgb255=(246,242,30) rgb255=(247,244,28) rgb255=(247,245,27) rgb255=(248,247,26) rgb255=(248,248,24) rgb255=(249,249,22) rgb255=(249,251,21) },
   colorbar horizontal,
   point meta min=0.000000e+00, point meta max=6.000000e+00,
   colorbar style={width=10cm, xtick={0.000000e+00,2.000000e+00,4.000000e+00,6.000000e+00}}
]
\addplot [draw=none] coordinates {(0,0)};
\end{axis}
\end{tikzpicture}  
\end{figure}
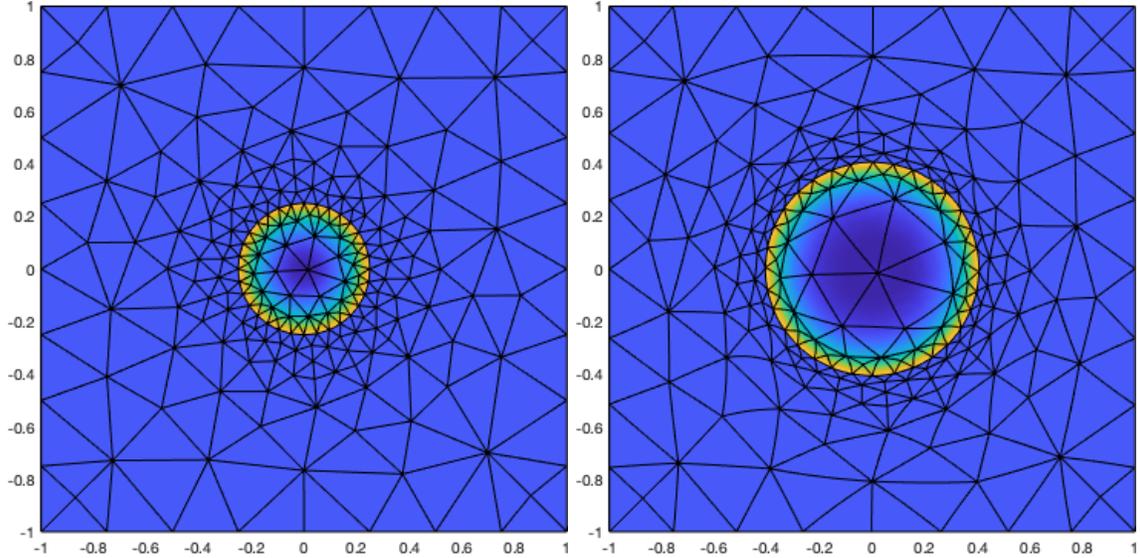

We take $\Omega = [-1, 1]^2$ as our two-dimensional spatial domain with the initial condition
given in terms of the density, velocity, and pressure by a function of 
the distance from the origin $r = \sqrt{x_1^2 + x_2^2}$ as

\begin{equation}
\begin{aligned}
 \rho(r) &= \begin{cases} \frac{5.378}{0.25^2}r^2 & r \leq 0.25 \\ 1 & \text{elsewhere} \end{cases} \\
 v(r) &= \begin{cases} \frac{1.304}{0.25}r & r  \leq 0.25 \\ 0 & \text{elsewhere} \end{cases} \\
 P(r) &= \begin{cases} \frac{0.978}{0.25^2}r^2 + 1& r  \leq 0.25 \\ 10^{-3} & \text{elsewhere} \end{cases}.
\end{aligned}
\end{equation}
Like the other numerical examples considered in this work, 
the initial condition is constructed such that the shock is exactly meshed.
This initial condition is inspired by the solution to the Sedov problem at a time 
when the primitive variables have already taken on a ``bowl-like" structure.
The shock tracking solution is computed using a DG discretization
on a mesh with $422$ simplex elements of degree $p=2$, $q=2$ and a DIRK$2$
temporal discretization with $140$ time steps and final time $T = 0.14$ (Figure~\ref{fig:blastwaveMOLsoln}).
In the shock tracking framework, the ability to track curved shocks
accurately while maintaining radial symmetry comes very naturally. 
Since the initial value for pressure in the ambient region is 
close to zero and much smaller than the value inside the shock,
a slight modification was made to the SQP solver to avoid negative 
values for energy in the ambient region. After each SQP step, 
the elements in the ambient region with a negative value for energy
are reinitialized with their element-wise average. 
This helps stabilize the solver in its initial iterations, 
and it quickly converges from there.

\section{Conclusion and future work}
\label{sec:conclude}

We extend the high-order implicit shock tracking (HOIST) 
framework for inviscid, steady conservation laws introduced in \cite{zahr2020implicit} 
to the unsteady case by a method of lines discretization via a 
diagonally implicit Runge-Kutta (DIRK) method. 
In essence, we are ``solving a steady problem at each timestep", 
and as such inherit the desirable qualities discussed 
for the steady case in our previous work, namely a conservative
and feature-aligned discretization at each timestep.
We demonstrate our framework on a variety of one- and two- dimensional
examples, and in particular show that our method is capable of preserving the 
design order of accuracy of our high-order temporal discretization
for both the solution and the shock location, the latter of which 
is inaccessible to shock capturing methods because they do not maintain
a perfect representation of discontinuities.
The ability to construct good initial guesses for the optimization problem
from information at the previous timestep is key to enable the SQP solver
to converge rapidly.
This is in contrast to the steady case, where the initial mesh generation
and corresponding guess for the solution has to be done independent of 
the shock location because we lack any \textit{a priori} knowledge.
Correspondingly, the steady case requires significant mesh deformation 
and edge collapses, which needs a more sophisticated line search procedure
based on the $\ell_1$ merit function to achieve convergence. 
The development of this method of lines approach gives 
us an additional tool to tackle unsteady problems in addition to the
previously developed space-time approach. Generally speaking,
the method of lines approach will be more practical and scale better
as the size and difficulty of the problem increases. However, it is limited 
in that it cannot handle colliding shocks (triple points in space-time) 
without complex mesh operations and solution reinitialization. 
In these cases, the space-time approach is preferred due to its generality 
of tracking discontinuities in space-time, which naturally handles triple points.

We see two important avenues of future work in the current setting of 
time-dependent, inviscid conservation laws. The first issue to consider
is shock formation, most commonly illustrated by the time 
evolution of Burgers' equation from a smooth initial condition or the
steepening of the waves trailing the shock in the Shu-Osher problem.
Implicit shock tracking has the potential to cleanly form and 
subsequently track shocks due to its $r$-adaptive behavior and the optimization
formulation that \textit{implicitly} tracks shocks and simultaneously resolves
the flow; this was shown by the authors of the MDG-ICE method for
inviscid and viscous shocks in the space-time setting
\cite{corrigan2019unsteady, kercher_moving_2021}.
In contrast, shock tracking approaches that handle shocks explicitly
do not readily show potential to handle shock formation without significant
specialization.
%It appears that other shock tracking approaches also currently avoid 
%dealing with shock formation, which occurs as a result of the steepening 
%of the trailing waves in the Shu-Osher problem 
%\cite{rawat2010high,corrigan2019unsteady}, just as we do here.
%However, we believe the HOIST framework has the potential to cleanly form and 
%subsequently track shocks due to its $r$-adaptive behavior and the optimization
%formulation that \textit{implicitly} tracks shocks and simultaneously resolves
%the flow. In contrast, shock tracking approaches that handle shocks explicitly
%do not readily show potential to handle shock formation without significant
%specialization.
%analogous to the modification of a standard flow solver to explicitly track the shock in the first place.
%To the best of our knowledge, no method has successfully demonstrated 
%the design order of accuracy for the full Shu-Osher problem or more complex flows featuring discontinuities.
%Our goal is to overcome the long accepted conventional 
%wisdom from the shock capturing community that only first order global accuracy
%is formally attainable in the presence of discontinuities, instead of the design order.

Another direction would be the incorporation of topology changes, such as 
edge flips, refinement and coarsening. In this work, we use a fixed mesh topology
as our shock moves across the domain. This limits our ability to 
accurately resolve the smooth flow in certain areas that only have a 
fixed amount of grid resolution to represent an increasingly larger
portion the spatial domain over time. The issue is somewhat exacerbated by the 
fact that we are working with very coarse meshes to begin with.
This difficultly manifests itself in our numerical experiments, such as 
in the area immediately trailing the shock of the Shu-Osher problem 
as well as inside the circular region of the blast wave example as the
shock expands radially; this loss of resolution is indicated by the
gradual increase in the objective function over time
(Figure~\ref{fig:shuosher_sqp}).
We emphasize that this is not a fundamental issue and the method
performs as expected under the constraint imposed by the fixed
topology. However, the incorporation of adaptive mesh coarsening
and refinement on curved meshes combined with $L^2$ error minimizing
solution transfer between meshes will unlock the full potential of the
method.
%to unlock the full potential of the method,
%the 
%framework works as expected and accurately as possible
%for the given fixed topology, but rather the addition of topology changes
%would help unlock the full potential of the method. We plan to develop an adaptive mesh 
%coarsening and refinement procedure on curved meshes so we can maintain
%the initial level of resolution over time. This also requires corresponding
%$L^2$ error minimizing solution transfer procedure to accompany
%each individual topological operation.

%Mention plan to further develop spacetime approach in tandem here?
Our implicit shock tracking framework has only been introduced and tested 
on inviscid conservation laws, where perfect discontinuities arise. 
We expect the framework to be useful for viscous problems with 
smooth, high-gradient solution features, where the second order 
viscous terms can be treated using standard techniques \cite{arnold2002unified}.
%There is already additional work in progress on the topic of 
%iterative solvers and preconditioners for the SQP solution of the linear system
%to make the approach practical for large-scale problems. 
Finally, this paper only considers relatively simple model problems 
in one- and two- dimensions which feature fairly tame shock motions.
Future work will deal with more complex shock phenomena, such
as shock-shock and shock-boundary layer interaction.
%We hope to demonstrate the usefulness of our proposed shock 
%tracking framework on more complex
%2D problems as well as relevant 3D problems in the near future.

\section*{Acknowledgments}
This material is based upon work supported by the Air Force Office of Scientific Research (AFOSR) 
under award numbers FA9550-20-1-0236, FA9550-22-1-0002, and FA9550-22-1-0004. 
The content of this publication does not necessarily reflect the position or 
policy of any of these supporters, and no official endorsement should be inferred. 

\bibliographystyle{plain}
\bibliography{biblio_shktrkuns}

\end{document}